\documentclass{amsart}
\usepackage{amsfonts}
\usepackage{amsmath,amssymb}
\usepackage{amsthm}
\usepackage{amscd}
\usepackage{graphics}
\usepackage{graphicx}
{
\theoremstyle{remark}
\newtheorem{Def}{{\rm Definition}}
\newtheorem{Ex}{{\rm Example}}
\newtheorem{Rem}{{\rm Remark}}

}
\newtheorem{Cor}{Corollary}
\newtheorem{Prop}{Proposition}
\newtheorem{Thm}{Theorem}

\begin{document}
\title[Differential topology of manifolds admitting round fold maps]{Round fold maps and the topologies and the differentiable structures of manifolds admitting explicit ones}
\author{Naoki Kitazawa}
\address{
}
\email{n-kitazawa@imi.kyushu-u.ac.jp}
\maketitle

\begin{abstract}
{\it {\rm (}Stable{\rm )} fold} maps are fundamental tools in a generalization of the theory of Morse functions on smooth manifolds
 and its application to studies of geometric properties of smooth manifolds. {\it Round} fold maps were introduced as stable fold
 maps such that the sets of all of the singular values of them are concentric spheres by the author
 in 2013-4. Topological properties of such maps and topological information of their source manifolds such as homology
 and homotopy groups have been studied under appropriate conditions by the author.

\ \ \ In this
 paper, we redefine round fold maps respecting the definition. As more precise information of manifolds admitting round fold maps, we study the topologies and differentiable structures of
 manifolds admitting such maps under appropriate differential topological conditions.   
\end{abstract}

\section{Introduction}
\label{sec:1}
\subsection{Historical backgrounds and what are presented in the present paper}
{\it Fold} maps are fundamental tools in considering the theory of higher dimensional versions of Morse functions on smooth manifolds and its application to
 studies of smooth manifolds. Studies of geometric (algebraic and differential topological) properties of fold maps and their source manifolds have been important. Such studies
 were started by Whitney (\cite{whitney}) and Thom (\cite{thom}) in the 1950s. A {\it fold} map from an $m$-dimensional closed smooth manifold into an $n$-dimensional smooth
 manifold without boundary ($m \geq n \geq 2$) is a smooth map
 such that at each singular point, the function is of the form
$$(x_1, \cdots, x_m) \mapsto (x_1,\cdots,x_{n-1},\sum_{k=n}^{m-i}{x_k}^2-\sum_{k=m-i+1}^{m}{x_k}^2)$$
 for an integer $0 \leq i \leq \frac{m-n+1}{2}$: the integer $i$ is uniquely
 determined and called the {\it index} of the singular point). A Morse function
 is naturally regarded as a fold map ($n=1$).
  For such a map, the following two hold.
\begin{enumerate}
\item The set consisting of all the singular points (the {\it singular set}) is a closed smooth submanifold of dimension $n-1$ of the source manifold. 
\item The restriction map to the singular set is a smooth immersion of codimension $1$.
\end{enumerate}
 
We also note that if the restriction map to the singular set is an immersion with normal crossings, then it is {\it stable}: {\it stable} maps are important in the theory of global singularity; see \cite{golubitskyguillemin} for example. 

Since around the 1990s, fold maps with additional conditions such
 as {\it special generic} maps, which were studied in \cite{burletderham}, \cite{furuyaporto}, \cite{saeki2}, \cite{saeki3}, \cite{saekisakuma} and \cite{sakuma}, and algebraic and differential topological properties of manifolds admitting such maps have been
 actively studied. A {\it special generic} map
 is defined as a fold map such that the indices of singular points are always $0$. A Morse function on a homotopy
 sphere with just
 two singular points is regarded as a special generic map; every smooth homotopy sphere of dimension $k \neq 4$ and the $4$-dimensional standard sphere $S^4$ admits such
 a function and a manifold admitting such a function is homeomorphic to a sphere (see \cite{milnor2} and \cite{milnor3} and see also \cite{reeb}). Moreover, we easily obtain a special generic map from any standard sphere of dimension $k_1 \geq 2$ into the $k_2$-dimensional Euclidean space ${\mathbb{R}}^{k_2}$ by a natural
 projection under the assumption that $k_1 \geq k_2 \geq 1$ holds. On the other hand, it was
 shown that a homotopy sphere of dimension $k_1$ admitting a special generic map
 into ${\mathbb{R}}^{k_2}$ is diffeomorphic to the standard sphere $S^{k_1}$ under the assumption that $1 \leq k_1-k_2 \leq 3$ holds (see \cite{saeki2} and Example \ref{ex:1} of this paper). In addition, in \cite{saeki2}
 and \cite{saekisakuma}, manifolds admitting special generic maps into ${\mathbb{R}}^2$ and ${\mathbb{R}}^3$ are completely
 or partially classified. From a Morse function and its singular points, as is well-known, we can know homology groups and some information
 on homotopy of the source manifold and from a special generic map and its singular points, we can
 know more precise information such as the topology and the differentiable structure of the source manifold.

Later, in \cite{kitazawa2}, {\it round} fold maps, which will be mainly studied in this paper, were introduced. A {\it round} fold map
 is defined as a fold
 map satisfying the following three.
\begin{enumerate}
\item The singular set is a disjoint union of standard spheres.
\item The restriction map to the singular set is an embedding.
\item The set consisting of all the singular values (the {\it singular value set}) of the map is a disjoint union of spheres embedded concentrically. 
\end{enumerate}

\begin{figure}
\includegraphics[width=50mm]{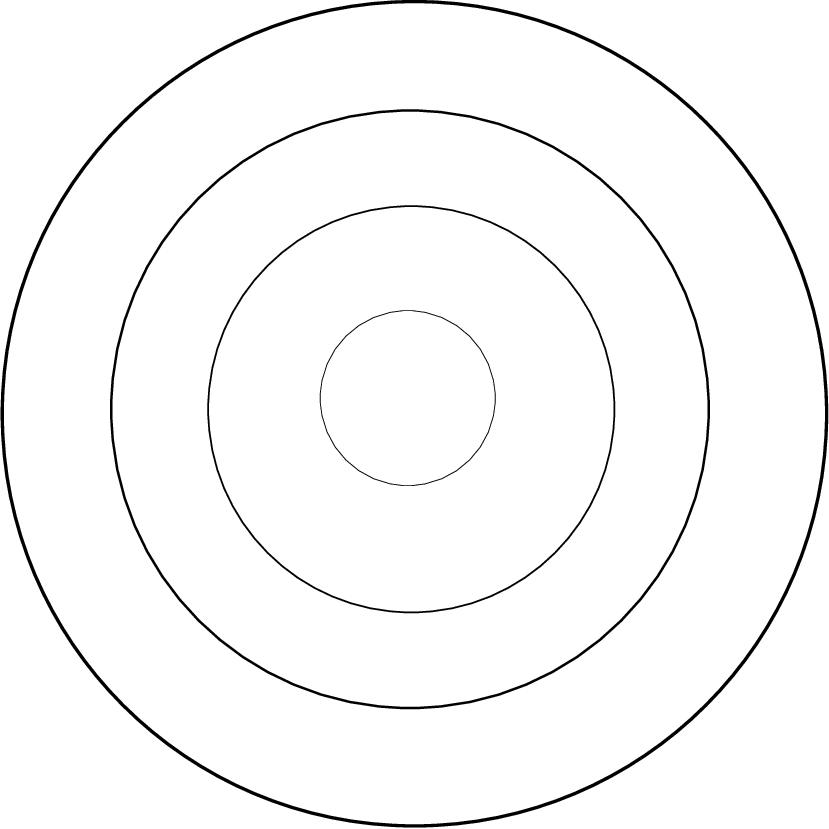}
\caption{The image (of the singular value set) of a round fold map into the plane.}
\label{fig:1}
\end{figure}
 For example, some special generic maps on homotopy spheres are
 round fold maps whose singular sets are connected. Any standard sphere whose dimension is $m>1$
 admits such a map into ${\mathbb{R}}^n$ under the assumption that the relation $m \geq n \geq 2$ holds and any homotopy sphere whose dimension is larger than $1$ and not $4$ admits
 such a map into the plane (see also \cite{saeki2} and Example \ref{ex:1} of this paper).

Homology and homotopy groups of manifolds admitting round fold maps were studied in \cite{kitazawa2} under appropriate
 conditions. For example, such groups of manifolds admitting round fold maps such that the inverse images of regular values are disjoint unions of
 spheres were studied. Such maps are generalizations of simplest special generic maps into Euclidean spaces whose singular sets are connected and standard spheres. 

 In this paper, as an advanced fundamental problem, we study the topologies and the differentiable structures of manifolds admitting round fold maps having appropriate differential topological structures.

\subsection{The content of the paper}
This paper is organized as the following.

In section \ref{sec:2}, first, we review the {\it Reeb space} of a smooth
 map, which is defined as the space consisting of all connected components of
 the inverse images of the map. Second, we reformulate {\it round} fold maps. For example, as a new ingredient, we define a function obtained by gluing two copies of a Morse function
 satisfying the following on the boundaries as a specific round fold map.  
\begin{itemize}
\item If the source manifold is a compact manifold with non-empty boundary, then the inverse image of the minimum is the boundary.  
\item Singular points are in the interior of the source manifold.
\item At distinct singular points, the values are distinct.
\end{itemize}

In section \ref{sec:3}, we first introduce known results on round fold maps on spheres and we construct round fold maps
 on $m$-dimensional closed smooth manifolds regarded as the total spaces of bundles over $S^n$ whose fibers are smooth manifolds and whose structure groups consist of diffeomorphisms such that the inverse images of axes are diffeomorphic
 to cylinders (Theorems \ref{thm:1}-\ref{thm:4}) under the assumption that the relation $m \geq n (\geq 2)$ holds. This is partially presented in \cite{kitazawa} and \cite{kitazawa2}.

In section \ref{sec:4}, we study the topologies and the differentiable structures of $m$-dimensional manifolds admitting round fold maps into ${\mathbb{R}}^n$ under the condition that the relation $m \geq 2n$ holds and additional algebraic and differential topological
 conditions. First, we recall
 a result on homology and homotopy groups of manifolds admitting round fold maps such that the inverse images of regular values are
 spheres shown in \cite{kitazawa} or \cite{kitazawa2} (Proposition \ref{prop:1}). 

 In subsection \ref{subsec:4.1}, under appropriate conditions, we construct
 a new round fold map from a manifold represented as a connected sum of two manifolds admitting round fold maps (Proposition \ref{prop:2}). Conversely, in
 subsection \ref{subsec:4.2}, we decompose a round fold map on a closed and connected manifold into two round fold maps so that a connected sum of
 the resulting source manifolds is the original source manifold (Proposition \ref{prop:4}). In these
 subsections, we apply generalized operations of surgery operations on stable maps from closed
 and simply-connected smooth manifolds into the plane which do not change the diffeomorphism types
 of the source manifolds in \cite{kobayashisaeki} ({\it R-operations}) and
 their inverse operations. Last, in subsection \ref{subsec:4.3}, by applying the obtained results and their proofs, we determine the topologies and the differentiable structures of manifolds admitting round fold
 maps such that the inverse images of their regular values are always disjoint unions of spheres mentioned before under a few additional
 differential topological conditions. More precisely, through Theorems \ref{thm:6}--\ref{thm:9}, we give a characterization of the family of closed and
 connected manifolds represented as connected sums of finite numbers of manifolds regarded as
 the total spaces of of bundles over $S^n$ whose fibers are diffeomorphic to $S^{m-n}$ and whose structure groups consist of diffeomorphisms by a family of round fold maps
 mentioned in the previous sentence. We also show several additional theorems stating that several fundamental closed and ($n-1$)-connected manifolds admit round fold maps into ${\mathbb{R}}^n$.
\subsection{Notation and terminologies}
 
 For two topological spaces $X \rightarrow Y$, subspaces $A \subset X$ and $B \subset Y$ and a homeomorphism $\phi:B \rightarrow A$, we denote the space obtained by attaching $Y$ to $X$ by $\phi$ by $X {\bigcup}_{\phi} Y$ and we omit the homeomorphism if the homeomorphism is clear in the context.  

 We also note on (homotopy) spheres. In this paper, an {\it almost-sphere} of dimension $k$ means a homotopy sphere
 given by gluing two $k$-dimensional standard
 closed discs together by a diffeomorphism between the boundaries. In other words, almost-sphere is a manifold represented as a {\it twisted double} of a standard closed disc. We can also define a {\it twisted double}
 of a compact and connected smooth manifold whose boundary is non-empty.

We often use terminologies on (fiber) bundles in this paper (see also \cite{steenrod}). For a topological space $X$, an {\it $X$-bundle} is a bundle
 whose fiber is $X$. A bundle whose structure group is $G$
 is said to be a {\it trivial} bundle if it is equivalent to the product bundle as a bundle whose structure group is $G$. Especially, a trivial bundle whose
 structure group is a subgroup of the homeomorphism group of the fiber is said to be a {\it topologically trivial} bundle. In this paper, a {\it smooth {\rm (}${\rm PL}${\rm )} bundle} means a bundle
 whose fiber is a smooth (resp. PL) manifold
 and whose structure group is a subgroup of the diffeomorphism group (resp. PL homeomorphism group) of the
 fiber. A {\it linear} bundle is
 a smooth bundle whose fiber is a standard disc (an unit disc) or a standard sphere (an unit sphere) and whose structure group consists of linear transformations on the fiber. 

 Throughout this paper, we assume that $M$ is a closed smooth manifold of dimension $m$, that
 $N$ is a smooth manifold of dimension $n$ without boundary, that $f:M \rightarrow N$ is a smooth map and that $m \geq n \geq 1$. Manifolds are smooth and of class $C^{\infty}$ and smooth maps between manifolds are also of class $C^{\infty}$ unless otherwise stated
 in the proceeding sections. 
 Note also that for a smooth map $c:X \rightarrow Y$, the {\it singular set} is, as explained, defined as the set of all singular points of $c$ and denoted by $s(c)$, the {\it singular value set} is defined as $c(s(c))$ and the {\it regular value set} is defined as $Y-c(s(c))$.

\subsection{Acknowledgement}
\thanks{This work is based on the doctoral dissertation by the author \cite{kitazawa2}, in which algebraic and differential topological
 properties of manifolds admitting round fold maps have been studied and the contents of the papers \cite{kitazawa} and \cite{kitazawa3} by the author are included. 

The author would like to express his gratitude to Mitsutaka Murayama, Osamu Saeki, Takahiro Yamamoto etc. for helpful comments and
constant encouragement.

The author is a member of the project Grant-in-Aid for Scientific Research (S) (17H06128 Principal Investigator: Osamu Saeki)
"Innovative research of geometric topology and singularities of differentiable mappings"

(https://kaken.nii.ac.jp/en/grant/KAKENHI-PROJECT-17H06128/)
 and supported by the project.
}

\section{Reeb spaces and round fold maps}
\label{sec:2}

\subsection{Reeb spaces} 

\begin{Def}
\label{def:1}
 Let $X$ and $Y$ be topological spaces. For $p_1, p_2 \in X$ and for a map $c:X \rightarrow Y$, 
 we define as $p_1 {\sim}_c p_2$ if and only if $p_1$ and $p_2$ are in
 a same connected component of $c^{-1}(p)$ for some $p \in Y$. ${\sim}_{c}$ is an equivalence relation.

We denote the quotient space $X/{\sim}_c$ by $W_c$. We call $W_c$ the {\it Reeb space} of $c$.
\end{Def}

 We denote the induced quotient map from $X$ into $W_c$ by $q_c$. We define $\bar{c}:W_c \rightarrow Y$
 so that the relation $c=\bar{c} \circ q_c$ holds.

$W_c$ is often homeomorphic to a polyhedron. For example, for a
 Morse function, the Reeb space
 is a graph and for a
 special generic map,
 the Reeb space is homeomorphic to a smooth manifold (see section 2 of \cite{saeki2}). 

In this section, we redefine {\it round} fold maps respecting the original definition (\cite{kitazawa}, \cite{kitazawa3} etc.). 

\subsection{Fundamental terms on round fold maps and constructing fundamental round fold maps}
\label{subsec:2.2}

First, we reformulate a {\it round} fold map.

Before that, we recall {\it $C^{\infty}$ equivalence} (see also \cite{golubitskyguillemin} for example). For two smooth maps $c_1:X_1 \rightarrow Y_1$ and $c_2:X_2 \rightarrow Y_2$, we say
 that they are {\it $C^{\infty}$ equivalent} or $c_1$ is said to be {\it $C^{\infty}$
 equivalent} to $c_2$ if there exist $C^{\infty}$ diffeomorphisms
 ${\phi}_X:X_1 \rightarrow X_2$ and ${\phi}_Y:Y_1 \rightarrow Y_2$ such that the following diagram commutes.

$$
\begin{CD}
X_1 @> {\phi}_X  >>X_2 \\
@VV c_1 V @VV c_2 V\\
Y_1 @> {\phi}_Y  >> Y_2
\end{CD}
$$

\begin{Def}[round fold map, \cite{kitazawa3}]
\label{def:2}
$f:M \rightarrow {\mathbb{R}}^n$ ($m \geq n \geq 2$) is said to be a {\it round} fold map if either of the follwoing hold.
\begin{enumerate}
\item $n=1$ holds and then $f$ is $C^{\infty}$ equivalent to
 a fold map $f_0:M_0 \rightarrow {\mathbb{R}}$ on a closed manifold $M_0$ such that the following three hold and we call the function a {\it normal form} of $f$.

\begin{enumerate}
\item $0$ is a regular value of $f_0$.  
\item Two Morse funtions defined as ${f_0} {\mid}_{{f_0}^{-1}(-\infty,0]}$ and ${f_0} {\mid}_{{f_0}^{-1}[0,+\infty)}$ are $C^{\infty}$ equivalent.
\item $f_0(S(f_0))$ is the set of all integers whose absolute values are positive and not larger than a positive integer.  
\end{enumerate}
\item $n \geq 2$ holds and $f$ is $C^{\infty}$ equivalent to
 a fold map $f_0:M_0 \rightarrow {\mathbb{R}}^n$ on a closed manifold $M_0$ such that the following three hold.

\begin{enumerate}
\item The singular set $S(f_0)$ is a disjoint union of standard spheres whose dimensions are $n-1$ and consists of $l \in \mathbb{N}$ connected components.
\item The restriction map $f_0 {\mid}_{S(f_0)}$ is an embedding.
\item Let ${D^n}_r:=\{(x_1,\cdots,x_n) \in {\mathbb{R}}^n \mid {\sum}_{k=1}^{n}{x_k}^2 \leq r \}$. Then, $f_0(S(f_0))={\sqcup}_{k=1}^{l} \partial {D^n}_k$ holds.  
\end{enumerate}

We call $f_0$ a {\it normal form} of $f$. We call a ray $L$ from $0 \in {\mathbb{R}}^n$ an {\it axis} of $f_0$ and
 ${D^n}_{\frac{1}{2}}$ the {\it proper core} of $f_0$. Suppose that for a round fold map $f$, its normal form $f_0$ and diffeomorphisms
 $\Phi:M \rightarrow M_0$ and $\phi:{\mathbb{R}}^n \rightarrow {\mathbb{R}}^n$, $\phi \circ f=f_0 \circ \Phi$ holds. Then
 for an axis $L$ of $f_0$, we also call ${\phi}^{-1}(L)$ an {\it axis} of $f$ and for the proper core ${D^n}_{\frac{1}{2}}$ of $f_0$, we
 also call ${\phi}^{-1}({D^n}_{\frac{1}{2}})$ a {\it proper core} of $f$.
 Note that in the case where $n=1$, we can similarly define a proper core and an axis of $f$ ($f_0$). 
\end{enumerate}
\end{Def}


In the case where the target space is $\mathbb{R}$ (and for example, where the source manifold is connected), a round fold map is regarded as a so-called {\it twisted double} of a Morse function. We can easily see that any closed and connected
 manifold represented as a twisted double of a compact and connected manifold admits such a function.

Let $f$ be a normal form of a round fold map into ${\mathbb{R}}^n$ ($n \geq 2$) and $P:={D^{n}}_{\frac{1}{2}}$. We set $E_1:=f^{-1}(P)$ and $E_2:=M-f^{-1}({\rm Int} P)$. We set
 $F:=f^{-1}(p)$
 for a point $p \in \partial P$. We put $Q:={\mathbb{R}}^n-{\rm Int} P$. Let $f_P:=f {\mid}_{E_1}:E_1 \rightarrow P$ if $F$ is non-empty and let $f_Q:=f {\mid}_{E_2}:E_2 \rightarrow Q$.

$f_P$ gives us a trivial bundle. ${f_P} {\mid}_{\partial E_1}:\partial E_1 \rightarrow \partial P$ gives
 us a trivial smooth bundle
 over $\partial P$ if $F$ is non-empty. $f_Q {\mid}_{\partial E_2}:\partial E_2 \rightarrow \partial Q$
 gives us a trivial smooth bundle over $\partial Q$. \\

We can give $E_2$ a bundle structure over $\partial Q$ as follows. \\ 

Since for ${\pi}_P(x):=\frac{1}{2} \frac{x}{|x|}$ ($x \in Q$), ${\pi}_P \circ f {\mid}_{E_2}$ is
 a proper submersion, this map makes $E_2$ a smooth $f^{-1}(L)$-bundle over $\partial Q$ (apply Ehresmann's fibration
 theorem \cite{ehresmann}): more precisely, in this situation, $L$ must be restricted to the set of all points such that the absolute values are not smaller than $\frac{1}{2}$. We call this bundle the {\it surrounding bundle} of $f$. Note that the structure group of this bundle
 is regarded as the group of diffeomorphisms on $f^{-1}(L)$ preserving the function $f {\mid}_{f^{-1}(L)}:f^{-1}(L) \rightarrow L (\subset \mathbb{R})$, which
 is naturally regarded as a Morse function.


For a round fold map $f$ into ${\mathbb{R}}^n$ ($n \geq 2$) which is not a normal form, we can consider similar maps and bundles. We call
 bundles naturally corresponding to the surrounding bundle of a normal form of $f$ a {\it surrounding bundle} of $f$.

\begin{figure}
\includegraphics[width=50mm]{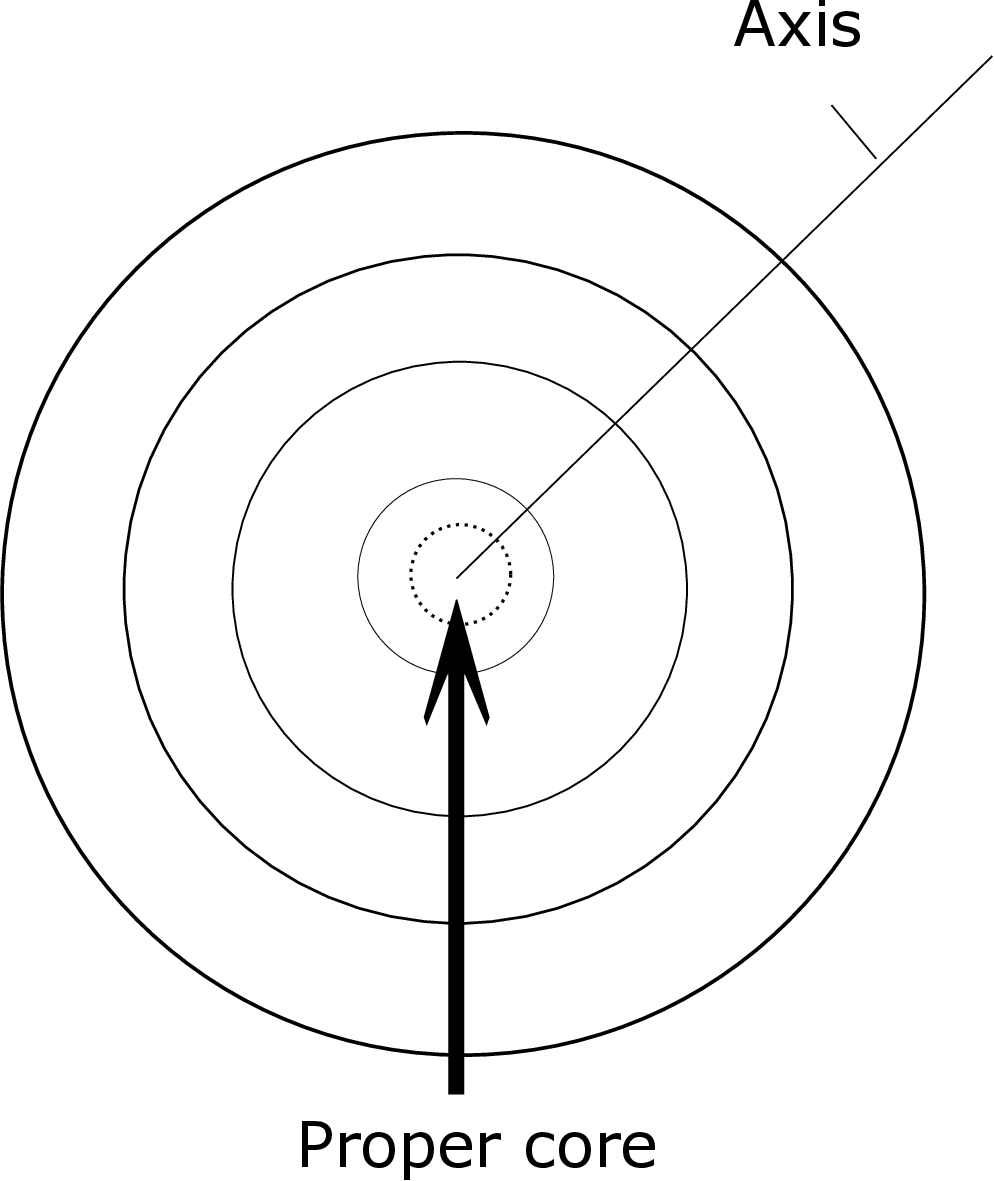}
\caption{A round fold map into the plane.}
\label{fig:2}
\end{figure}

%
%
%
For these terminologies, see FIGURE \ref{fig:2}.
We can define the following condition for a round fold map. 

\begin{Def}
\label{def:3}
For an integer $n \geq 2$, let $f:M \rightarrow {\mathbb{R}}^n$ be a round fold map. 
If a surrounding bundle of $f$ as above 
 is a topologically trivial bundle, then
 $f$ is said to be {\it topologically trivial}. If the bundle is a trivial PL bundle, then $f$ is said to be {\it PL trivial}. If the bundle is a trivial smooth bundle, then $f$ is said to
 be {\it $C^{\infty}$ trivial}.  
\end{Def}

We can construct a round fold map into $\mathbb{R}$ or a $C^{\infty}$ trivial round fold map into ${\mathbb{R}}^n$ ($n \neq 1$) as in the
following manner.

Let $\bar{M}$ be a compact manifold with non-empty boundary $\partial \bar{M}$. Then, there
 exists a Morse
 function $\tilde{f}:\bar{M} \rightarrow [a,+\infty)$ such that the following three hold.

\begin{enumerate}
\item $a$ is the minimum of $\tilde{f}$. 
\item ${\tilde{f}}^{-1}(a)= \partial \bar{M}$.
\item Singular points of $\tilde{f}$ are always in the interior ${\rm Int} \bar{M}$ and at distinct singular points, the values are always distinct.
\end{enumerate}

Let $\Phi:\partial (\bar{M} \times \partial ({\mathbb{R}}^n-{\rm Int} D^n)) \rightarrow \partial (\partial \bar{M} \times D^n)$
 and $\phi:\partial ({\mathbb{R}}^n-{\rm Int} D^n) \rightarrow \partial D^n$ be diffeomorphisms. Let $p_1:\partial \bar{M} \times \partial ({\mathbb{R}}^n-{\rm Int} D^n) \rightarrow \partial ({\mathbb{R}}^n-{\rm Int} D^n)$ and
 $p_2:\partial \bar{M} \times \partial D^n \rightarrow \partial D^n$ be
 the canonical projections. Suppose that the following diagram commutes. 

$$
\begin{CD}
\partial \bar{M} \times \partial ({\mathbb{R}}^n-{\rm Int} D^n)  @> \Phi >> \partial \bar{M} \times \partial D^n \\
@VV p_1 V @VV p_2 V \\
\partial ({\mathbb{R}}^n-{\rm Int} D^n) @> \phi >> \partial D^n
\end{CD}
$$

By using the diffeomorphism $\Phi$, we construct a manifold $M:=(\partial \bar{M} \times D^n) {\bigcup}_{\Phi} (\bar{M} \times \partial ({\mathbb{R}}^n-{\rm Int} D^n))$. By using the
 diffeomorphism $\phi$, we can construct $D^n {\bigcup}_{\phi} ({\mathbb{R}}^n-{\rm Int} D^n)$, which is diffeomorphic to ${\mathbb{R}}^n$. 
Let $p:\partial \bar{M} \times D^n \rightarrow D^n$ be the canonical projection. Then, by gluing the two maps $p$ and $\tilde{f} \times {\rm id}_{\partial D^n}$
 together by the pair of diffeomorphisms $(\Phi,\phi)$, a round fold map $f:M \rightarrow {\mathbb{R}}^n$ is obtained; in this
 situation, we regard ${\mathbb{R}}^n-{\rm Int} D^n$ as $\partial D^n \times [a,+\infty)$ by a diffeomorphism between the spaces.

If $\bar{M}$ is a compact manifold without boundary, then there exists a Morse
 function $\tilde{f}:\bar{M} \rightarrow [a,+\infty)$ such that $\tilde{f}(\bar{M}) \subset (a,+\infty)$ and that at distinct singular points, the values are distinct. We are enough
 to consider $\tilde{f} \times {\rm id}_{S^{n-1}}$
 and embed $[a,+\infty) \times S^{n-1}$ into ${\mathbb{R}}^n$ to construct a round fold map whose source manifold is $\bar{M} \times S^{n-1}$.

We call the construction of a round fold map $f$ here a {\it trivial spinning construction}.

\section{Round fold maps on spheres and bundles over standard spheres}
\label{sec:3}
In this section, we introduce fundamental discussions on round fold
 maps on homotopy spheres which follow directly from known results and characterize manifolds regarded as total spaces of bundles over standard spheres by
 certain round fold maps.

The following example is from fundamental discussions of \cite{saeki2}. 

\begin{Ex}
\label{ex:1}
\begin {enumerate}
\item
\label{ex:1.1}
 Let $M$ be a closed manifold of dimension $m$. Let $n$ be an integer such that $m> n \geq 2$ and $n \neq 4, 5$ hold. A round fold map $f:M \rightarrow {\mathbb{R}}^n$ whose singular set is connected exists if and only if $M$ is a homotopy sphere
 admitting a special generic map into ${\mathbb{R}}^n$ whose Reeb space is homeomorphic to $D^n$. 
\item
\label{ex:1.2}
 Any standard sphere of dimension $m$ admits a map into ${\mathbb{R}}^n$ as above if $m \geq n \geq 2$ holds. Furthermore, any homotopy sphere of dimension $m>1$ admits a map into the plane
 as above unless $m=4$ according to a discussion in section 5 of \cite{saeki2}.    
\item
\label{ex:1.3}
 Let $m$ be an integer larger than $3$ and $n$ be an integer satisfying $m-n=1,2,3$. In section 4 of \cite{saeki2} and \cite{saeki3}, it is shown that
 if a homotopy sphere of dimension $m$ admits a special generic map into ${\mathbb{R}}^n$, then the homotopy sphere is a standard sphere.

Thus, in the situation of this part, if on a homotopy sphere of dimension $m$, a round fold map into ${\mathbb{R}}^n$ whose singular set is connected exists, then the sphere is diffeomorphic to $S^m$.

There are $28$ types of  $7$-dimensional oriented homotopy spheres modulo orientation preserving diffeomorphism, which we can know from \cite{eellskuiper} and \cite{milnor}, and recently, Wrazidlo has shown that $7$-dimensional oriented homotopy spheres of $14$ types do not admit special generic maps into ${\mathbb{R}}^3$ (\cite{wrazidlo}). 
\end{enumerate} 
\end{Ex} 

Note also that round fold maps in Example \ref{ex:1} are topologically trivial. They are also
 PL trivial and $C^{\infty}$ trivial (see also Example 3 (1) of \cite{kitazawa2}). 

The following theorem has been partially proven in \cite{kitazawa} as introduced in Example \ref{ex:2} (\ref{ex:2.1}) later.

\begin{Thm}
\label{thm:1}
Let $M$ be a closed manifold of dimension $m$. Let $n$ be an integer satisfying the relation $m \geq n \geq 2$. 
\begin{enumerate}
\item
\label{thm:1.1}
 Let $M$ be the total space of a smooth bundle over $S^n$ whose fiber is a closed manifold $F$. Then, there exists
 a $C^{\infty}$ trivial round fold map $f:M \rightarrow {\mathbb{R}}^n$ such that the inverse image of a point in a proper core of $f$ is diffeomorphic
 to a disjoint union of two copies of $F$, which are regarded as fibers of the $F$-bundle over $S^n$, and that $f^{-1}(L)$ is diffeomorphic to $F \times [0,1]$ for an axis $L$ of $f$. 
\item
\label{thm:1.2}
 Suppose that a topologically trivial
 round fold map $f:M \rightarrow {\mathbb{R}}^n$ exists and that for an axis $L$ of $f$ and a closed manifold $F$ of dimension $m-n$, $f^{-1}(L)$ is diffeomorphic
 to $F \times [0,1]$. Then, $M$ is the total space of an $F$-bundle over $S^n$. If $f$ is normally PL trivial, then $M$ is the total space of a PL $F$-bundle over $S^n$ and if $f$ is normally $C^{\infty}$ trivial, then $M$ is the total space of a smooth $F$-bundle over $S^n$.
\end{enumerate}
\end{Thm}

\begin{proof}
 We prove the first part.

 We may represent $S^n$ as $(D^n \sqcup D^n) \bigcup (S^{n-1} \times [0,1])$, where
 we identify $\partial (D^n \sqcup D^n)=S^{n-1} \sqcup S^{n-1}$ and $\partial (S^{n-1} \times [0,1])=S^{n-1} \sqcup S^{n-1}$. For a diffeomorphism $\Phi$ from $S^{n-1} \times (F \sqcup F)$ onto $\partial D^n \times (F \sqcup F)$ which is a bundle isomorphism between
 the trivial $F$-bundles over $\partial (D^n \sqcup D^n)=\partial (S^{n-1} \times [0,1])= S^{n-1} \sqcup S^{n-1}$ inducing the identification between
 the base spaces, we may represent $M$ as $((D^n \sqcup D^n) \times F) {\bigcup}_{\Phi} (S^{n-1} \times [0,1] \times F)=(D^n \times (F \sqcup F)) {\bigcup}_{\Phi} (S^{n-1} \times [0,1] \times F)$.

 There exists a Morse function $\tilde{f}:F \times [0,1] \rightarrow [a,+\infty)$, where $a \in \mathbb{R}$ is the minimal value.
\begin{enumerate}
\item $a$ is the minimum of $\tilde{f}$. 
\item ${\tilde{f}}^{-1}(a)=F \times \{0,1\}$.
\item Singular points of $\tilde{f}$ are always in the interior $F \times (0,1) \subset F \times [0,1]$ and at distinct singular points, the values are always distinct.
\end{enumerate}

We consider a map $\tilde{f} \times {\rm id}_{S^{n-1}}$ and the canonical projection $p:D^n \times (F \sqcup F) \rightarrow D^n$. For the maps $\Phi$, $\tilde{f} \times {\rm id}_{S^{n-1}}$, $p$
 and a diffeomorphism $\phi:\partial ({\mathbb{R}}^n-{\rm Int} D^n) \rightarrow \partial D^n$, the following diagram commutes.

$$
\begin{CD}
F \times (\{0\} \sqcup \{1\}) \times \partial ({\mathbb{R}}^n-{\rm Int} D^n) @> \Phi >> F \times (\{0\} \sqcup \{1\}) \times \partial D^n \\
@VV {\tilde{f}} {\mid}_{F \times (\{0\} \sqcup \{1\})} \times {\rm id}_{\partial ({\mathbb{R}}^n-{\rm Int} D^n)} V @VV p {\mid}_{F \times (\{0\} \sqcup \{1\}) \times \partial D^n} V \\
\{a\} \times \partial ({\mathbb{R}}^n-{\rm Int} D^n) @> \phi >> \partial D^n  
\end{CD}
$$  

Then, by gluing the two maps $p$ and $\tilde{f} \times {\rm id}_{\partial ({\mathbb{R}}^n-{\rm Int} D^n)}$ together by using diffeomorphisms
 $\Phi$ and $\phi$ or using a trivial spinning construction introduced just after Definition \ref{def:3}, we have a $C^{\infty}$ trivial round fold map $f:M \rightarrow {\mathbb{R}}^n$. We see that $f$ is a round fold map satisfying the given
 conditions.

 We now prove the second part.

 Suppose that there exists a topologically trivial round fold map $f:M \rightarrow {\mathbb{R}}^n$ and that for an axis $L$ of $f$ and a closed
 manifold $F$, $f^{-1}(L)$ is diffeomorphic to $F \times [0,1]$. Then, for a
 diffeomorphism $\Phi$ from $S^{n-1} \times (F \sqcup F)$ onto $\partial D^n \times (F \sqcup F)$ which is a bundle isomorphism between
 the trivial ($F \sqcup F$)-bundles over $S^{n-1}=\partial D^n$ inducing a diffeomorphism between the base spaces, $M$ is regarded
 as $(D^n \times (F \sqcup F)) {\bigcup}_{\Phi} (S^{n-1} \times [0,1] \times F)$ in the topology category. Thus, $M$ is the
 total space of an $F$-bundle over $S^n$. 

 Moreover, suppose that there exists a PL ($C^{\infty}$) trivial round fold map $f:M \rightarrow {\mathbb{R}}^n$ and
 that for an axis $L$ of $f$ and a closed
 manifold $F$, $f^{-1}(L)$ is diffeomorphic to $F \times [0,1]$. In this case, in the PL (resp. $C^{\infty}$)
 category, for a
 diffeomorphism $\Phi$ from $S^{n-1} \times (F \sqcup F)$ onto $\partial D^n \times (F \sqcup F)$ which is a bundle isomorphism between
 the trivial ($F \sqcup F$)-bundles over $S^{n-1}=\partial D^n$ inducing a diffeomorphism between
 the base spaces, $M$ is regarded
 as $(D^n \times (F \sqcup F)) {\bigcup}_{\Phi} (S^{n-1} \times [0,1] \times F)$.

 This completes the proof of both parts of the theorem.
\end{proof}
\begin{Rem}
For a closed and connected manifold $F$, twisted doubles of the cylinder $F \times [-1,1]$ are characterized as the total spaces of $F$-bundles over $S^1$. We can easily reach the version where the target space is $\mathbb{R}$ of Theorem \ref{thm:1}.
\end{Rem}

\begin{Ex}
\label{ex:2}
\begin{enumerate}
\item
\label{ex:2.1} 
Let $M$ be a closed manifold of dimension $m$ and let $M$ be the total space of a smooth bundle over $S^n$ whose
 fiber is an almost-sphere $\Sigma$ of dimension $m-n$ satisfying the relation $m>n \geq 2$. Then, there exists
 a $C^{\infty}$ trivial round fold map $f:M \rightarrow {\mathbb{R}}^n$ such that the inverse image of a point in a proper core of $f$ is diffeomorphic
 to a disjoint union of two copies of $\Sigma$ and that $S(f)$ consists of $2$ connected components and is the disjoint union
 of the set of
 all fold points of index $0$ and the set of
 all fold points of index $1$. This explanation is also presented in \cite{kitazawa} and the maps here are important ones in this paper. See FIGURE \ref{fig:3}, in which such a map in the $n=1$ case and one in the $n \neq 1$ case are presented with their Reeb spaces.
\begin{figure}
\includegraphics[width=50mm]{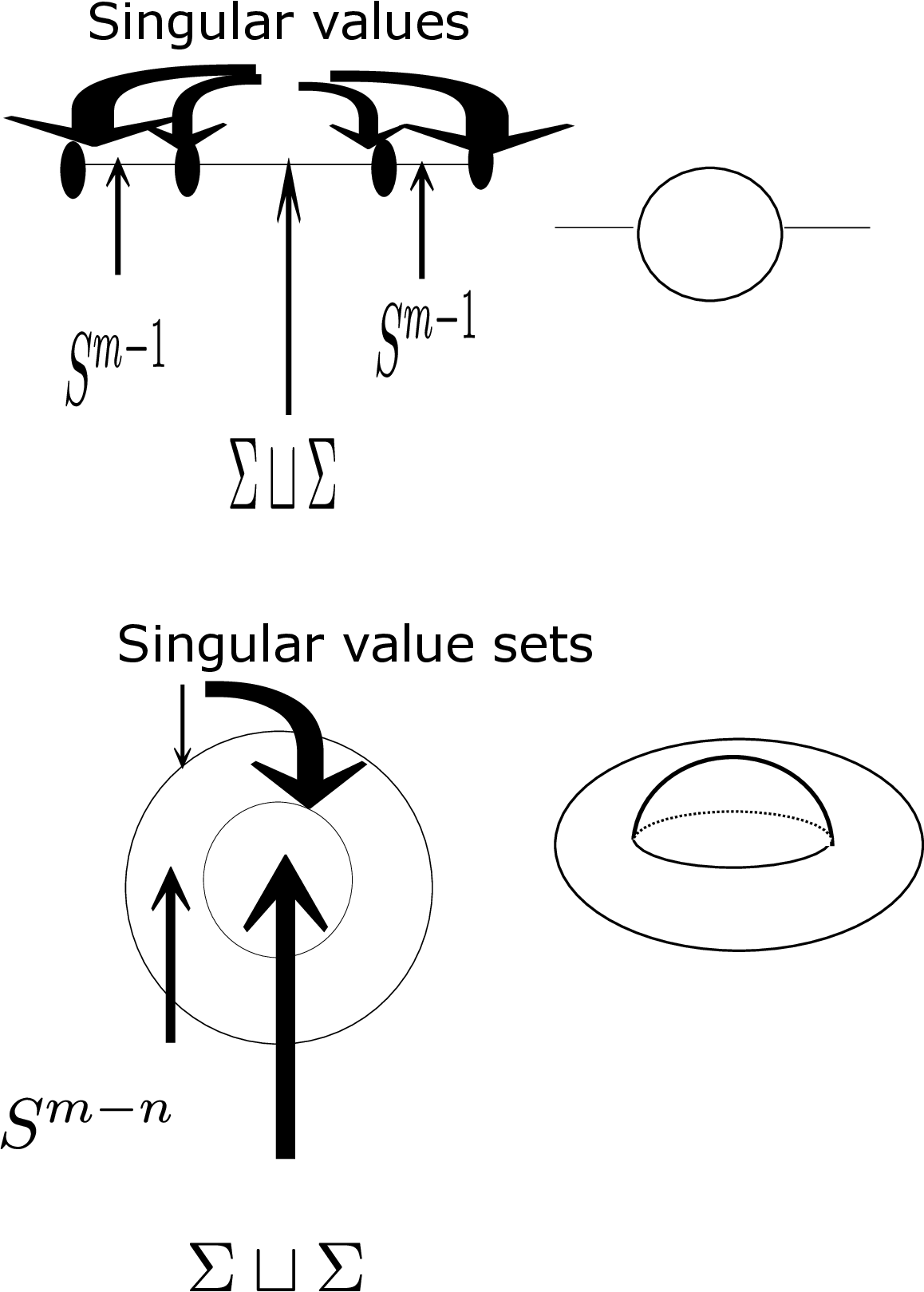}
\caption{Explicit round fold maps of Example \ref{ex:2} (\ref{ex:2.1}) (inverse images of regular values are as noted) and their Reeb spaces.}
\label{fig:3}
\end{figure}

\item
\label{ex:2.2}
 Let $k \in \mathbb{N}$ and $k \geq 3$. Let $M$ be the orthogonal group $O(k)$ of dimension $k$ or the special orthogonal group $SO(k)$ of dimension $k$. By considering
 a natural inclusion from the orthogonal group $O(k-1)$ of dimension $k-1$ into $O(k)$ or a natural inclusion from the rotation group $SO(k-1)$ of dimension $k-1$ into $SO(k)$,
 we obtain the natural homogeneous manifold diffeomorphic to the standard sphere $S^{k-1}$. Thus, $M$ is regarded as the total space of a smooth $O(k-1)$-bundle over $S^{k-1}$ if it
 is $O(k)$ and it is regarded as the total space of a smooth $SO(k-1)$-bundle if it is $SO(k)$. We obtain a round fold map from $M$ into ${\mathbb{R}}^{k-1}$ as
 presented in Theorem \ref{thm:1}.  
\end{enumerate}
\end{Ex}

 

We have the following theorem by virtue of the fact that there exist homotopy spheres homeomorphic
 but not diffeomorphic to $S^7$ and regarded as total spaces of linear $S^3$-bundles over
 $S^4$ (see \cite{eellskuiper} and \cite{milnor}) and Example \ref{ex:2}, which is presented in \cite{kitazawa}.

\begin{Thm}[\cite{kitazawa}]
\label{thm:2}
There exist homotopy spheres that are homeomorphic but not diffeomorphic to $S^7$ and they admit
 $C^{\infty}$ trivial round fold maps into ${\mathbb{R}}^4$ whose singular sets consist of $2$ connected components. The $7$-dimensional
 standard sphere $S^7$, which is the total space of a linear $S^3$-bundle over $S^4$, also admits such a map. 
\end{Thm}

\begin{Rem}
\label{rem:2}
If there exists a special generic map from a homotopy sphere homeomorphic to $S^7$
 into ${\mathbb{R}}^4$, then the homotopy sphere is diffeomorphic to $S^7$ from Example \ref{ex:1} (\ref{ex:1.3}). This
 means that we cannot reduce the numbers of connected components of the singular sets of
 round fold maps in Theorem \ref{thm:2} for homotopy spheres not diffeomorphic to $S^7$. 
\end{Rem}

In the case where $n=2$ is assumed, we have the following theorems.

\begin{Thm}
\label{thm:3}
Let $M$ be a closed manifold of dimension $m \geq 7$, $m=3$ or $m=4$. If $m \geq 7$ is assumed, let $F$ be a
 closed and simply-connected manifold of dimension $m-2$ and if $m=3,4$ is assumed, let $F$ be the {\rm (}$m-2${\rm )}-dimensional standard sphere $S^{m-2}$, then, the following two are equivalent.
\begin{enumerate}
\item \label{thm:3.1}
 $M$ is the total space of a smooth $F$-bundle over $S^2$.
\item \label{thm:3.2}
 $M$ admits a round fold map $f:M \rightarrow {\mathbb{R}}^2$ satisfying the following two conditions.
\begin{enumerate}
\item The inverse image of a point in a proper core of $f$ is diffeomorphic to a disjoint union of two copies of $F$.
\item For an axis $L$ of $f$, $f^{-1}(L)$ is diffeomorphic to $F \times [0,1]$.
\end{enumerate}
\end{enumerate} 
\end{Thm}
\begin{proof}
If $M$ admits a round fold map $f:M \rightarrow {\mathbb{R}}^2$ as in the condition (\ref{thm:3.2}), then since for a proper core $P$ of
 $f$, ${f} {\mid}_{f^{-1}(P)}: f^{-1}(P) \rightarrow P$ makes $f^{-1}(P)$ a trivial
 bundle, $f$ is $C^{\infty}$
 trivial by virtue of the pseudoisotopy theorem \cite{cerf} ($m \geq 7$) or the fact that a smooth ($S^{m-2} \times [0,1]$)-bundle over $S^1$ whose
 subbundle with the fiber $S^{m-2} \times \{0,1\} \subset S^{m-2} \times [0,1]$ is trivial is also trivial ($m=4$).

From Theorem \ref{thm:1}, this completes the proof. 
\end{proof}

We also have the following theorem, which is a generalization of Theorem \ref{thm:1} (\ref{thm:1.1}).  

\begin{Thm}
\label{thm:4}
Let $n$ be a positive integer.
Let $M$ be a manifold obtained by the following steps.

\begin{enumerate}
\item Let $F_1$ be a closed manifold and let $M_1$ be a closed manifold regarded as the total space of a smooth $F_1$-bundle over $S^n$.
\item Let $F_2$ be a closed manifold. Let $M_2$ be a closed manifold regarded as the total space of a smooth $F_2$-bundle over $M_1$ such that the restriction to any fiber of the previous
 bundle $M_1$ is a trivial smooth bundle.
\item For $M_{k}$ {\rm (}$k \geq 2${\rm )}, let $F_{k+1} \neq \emptyset$ be a closed manifold and we define a closed manifold $M_{k+1}$ as a manifold
 regarded as the total space of a smooth $F_{k+1}$-bundle over $M_k$ such that the restriction of the bundle to the product of every fiber
 appearing in the family $\{F_j\}_{j=1}^{k}$ of fibers of bundles $\{M_j\}_{j=1}^{k}$ is a trivial smooth bundle. We perform the methods inductively to
 obtain a manifold $M_l$ {\rm (}$l \geq 2${\rm )}.      
\item Let $M:=M_l$ {\rm (}$l \geq 2${\rm )}.
\end{enumerate}
In this situation, $M$ admits a $C^{\infty}$ trivial round fold map $f:M \rightarrow {\mathbb{R}}^n$ in the case $n \neq 1$ is assumed such that the inverse image
 of its axis is diffeomorphic to the product of all the $k$ manifolds belonging to the family $\{F_j\}_{j=1}^{k}$ and the closed interval $[-1,1]$. In addition, $M$ admits a round fold map $f:M \rightarrow {\mathbb{R}}^n$ represented as a twisted double of a Morse function on the product of all the $k$ manifolds belonging to the family $\{F_j\}_{j=1}^{k}$ and the closed interval $[-1,1]$ in the case $n =1$ is assumed.
\end{Thm}
\begin{proof}
In the case where $M=M_1$ holds, the result follows from Theorem \ref{thm:1} (\ref{thm:1.1}).

 Assume that $M_k$ admits a round fold map $f_k:M_k \rightarrow {\mathbb{R}}^n$ mentioned in the
 statement. Note that by Theorem \ref{thm:1} (\ref{thm:1.2}), $M_k$ is regarded as the total space of a smooth bundle over $S^n$ whose fiber is diffeomorphic
 to the product of all the $k$ manifolds belonging to the family $\{F_j\}_{j=1}^{k}$.

 We consider $M=M_{k+1}$ regarded as the total space of a
 smooth $F_{k+1}$-bundle over $M_k$ as in the mentioned steps. If we restrict the $F_{k+1}$-bundle to the inverse image ${f_k}^{-1}(P)$ of a proper core $P$ of $f_k$, which is regarded as the total space
 of a smooth trivial bundle over $P$ with a fiber diffeomorphic
 to the product of all the $k$ manifolds belonging to the family $\{F_j\}_{j=1}^{k}$, then
 it is a trivial smooth bundle by a condition on the structure of the bundle $M_{k+1}$ over $M_k$. If we restrict the $F_{k+1}$-bundle to the inverse image $f^{-1}({\mathbb{R}}^n-{\rm Int} P)$, which is
 diffeomorphic to the product of $\partial P$ or $S^{n-1}$, all the $k$ manifolds belonging to the family $\{F_j\}_{j=1}^{k}$ and the
 closed interval $[-1,1]$, then the resulting bundle is trivial by the same condition on the structure of the bundle $M_{k+1}$ over $M_k$.  

 From this discussion, $M_{k+1}$ is regarded as the total space of a smooth bundle over $S^n$ whose fiber is diffeomorphic
 to the product of all the $k+1$ manifolds belonging to the family $\{F_j\}_{j=1}^{k+1}$. Thus, we
 can construct a desired round fold map $f_{k+1}:M_{k+1} \rightarrow {\mathbb{R}}^n$ by using a method of the proof
 of Theorem \ref{thm:1} (\ref{thm:1.1}).

 By the induction, this completes the proof.     
\end{proof}

\begin{Ex}
\label{ex:3}
\begin{enumerate}
\item
\label{ex:3.1}
As discussed in Example \ref{ex:2} (\ref{ex:2.1}), we consider a natural inclusion from the special orthogonal group $SO(k_1)$ of dimension $k_1$
 into the special orthogonal group $SO(k_2)$ of dimension $k_2$ where $1 \leq k_1<k_2$ holds. We obtain a natural homogeneous
 manifold $SO(k_2)/SO(k_1)$. 

  Let $k_1,k_2$ and $k_3$ be integers satisfying $1 \leq k_1<k_2<k_3$. Then, by the canonical projection, the homogeneous manifold $SO(k_3)/SO(k_1)$ is regarded as a smooth bundle
 over the homogeneous manifold\\
 $SO(k_3)/SO(k_2)$ whose fiber is the homogeneous manifold\\ $SO(k_2)/SO(k_1)$ (see Part I 7 of \cite{steenrod} for such a bundle).
 
Moreover, assume that the bundle $SO(k_3)=SO(k_3)/SO(1)$

 over the homogeneous manifold $SO(k_3)/SO(k_2)$ is trivial and as a result regarded as $SO(k_2) \times SO(k_3)/SO(k_2)$. In this case, by setting $F_1=SO(k_3)/SO(k_2)$, $F_2=SO(k_2)/SO(k_1)$ and $n=k_3$, we may
 apply Theorem \ref{thm:4} and as a result obtain a round fold map from $SO(k_3+1)/SO(k_1)$ into ${\mathbb{R}}^{k_3}$. 

In this
 situation, let $k_1=1,2$, $k_2=3$ and $k_3=4$. It follows that the homogeneous manifold $SO(5)/SO(k_1)$ is regarded as a smooth bundle
 over the homogeneous manifold $SO(5)/SO(3)$ whose fiber is the homogeneous manifold $SO(3)/SO(k_1)$. The bundle $SO(4)=SO(4)/SO(1)$ is a trivial smooth bundle over the homogeneous
 manifold $SO(4)/SO(3)$, which is diffeomorphic to $S^3$. Thus, in the situation of Theorem \ref{thm:4}, by setting $F_1=SO(4)/SO(3)$. $F_2=SO(3)/SO(k_1)$ and $n=k_3=4$, we may
 apply Theorem \ref{thm:4} and as a result obtain a round fold map from $SO(5)/SO(k_1)$ into ${\mathbb{R}}^4$.  

\item
\label{ex:3.2}
 Let $M=S^2$. By considering the Whitney sum of a complex line bundle and a trivial complex line bundle over $M$, we obtain a complex vector bundle of dimension $2$ and by considering a natural quotient space of the vector
 space of the fiber, we have a smooth bundle over $S^2$ whose fiber $F_1$ is diffeomorphic to the
 complex projective space $\mathbb{C}P^1$ or $S^2$ (we consider the {\it projectivization} of the bundle). We denote the total space by $M_1$
 and similarly we consider the Whitney sum of a complex line bundle and a trivial complex line bundle over $M_1$ and its projectivization, whose fiber $F_2$ is also diffeomorphic to $S^2$. We denote the resulting total space, which
 is a $6$-dimensional connected and closed manifold, by $M_2$. Note that $M_1$ is a {\it 2-stage Bott manifold}
 and $M_2$ is a {\it 3-stage Bott manifold}. For {\it Bott manifolds}, see \cite{choimasudasuh} for example. We comment on
 the bundle structures and the diffeomorphism types of the manifolds $M_1$ and $M_2$ according to explanations of \cite{choimasudasuh}. In this comment, we use terminologies on the theory of characteristic classes
 of vector bundles such as {\it Chern classes} and {\it Pontrjagin classes} and fundamental theory of these characteristic classes (see \cite{milnorstasheff} for example).

  $M_1$ is regarded as the total space of a smooth $S^2$-bundle over $M=S^2$. There are just two isomorphism classes of smooth $S^2$-bundles over $M=S^2$ and the corresponding total spaces are not diffeomorphic mutually. The homology groups of these manifolds coincide; these manifolds
 are simply-connected and we have $H_2(M_1;\mathbb{Z}) \cong H_2(S^2 \times S^2;\mathbb{Z}) \cong \mathbb{Z} \oplus \mathbb{Z}$. Let $\nu$ be a generator of the cohomology group $H^2(M;\mathbb{Z}) \cong \mathbb{Z}$. More precisely, on
 each cohomology class $a\nu \in H^2(M;\mathbb{Z})$, we
 can consider the Whitney sum of a complex line bundle over $M$ whose Chern class is $a\nu \in H^2(M;\mathbb{Z})$ and a trivial complex line bundle over $M$ and its projectivization, whose fiber $F_1$ is diffeomorphic to $S^2$; if $a$ is even, then the resulting bundle is trivial and if not, then the bundle is not trivial. The
 cohomology group $H^2(M_1;\mathbb{Z})$ is isomorphic to $\mathbb{Z} \oplus \mathbb{Z}$ and freely generated by two elements $\alpha$ and $\beta$ satisfying the following; $\alpha$ is the element defined as
 the pull-back of a generator of the cohomology group $H^2(S^2;\mathbb{Z})$ of the base space of the bundle $M_1$ by the projection of the bundle $M_1$ and $\beta$ is
 an element such that the following two hold.
\begin{enumerate}
\item The value on the cycle representing a generator of the homology group $H_2(S^2;\mathbb{Z})$ of a fiber of the $S^2$-bundle $M_1$ is $1$.
\item The value is $0$ on any cycle representing an element of the group $H_2(M_1;\mathbb{Z}) \cong \mathbb{Z} \oplus \mathbb{Z}$ such that
 this cycle and the previous cycle generate the group $H_2(M_1;\mathbb{Z}) \cong \mathbb{Z} \oplus \mathbb{Z}$.  
\end{enumerate}
 To obtain $M_2$, we consider the Whitney sum of
 a complex line bundle and a trivial complex line bundle over $M_1$. On each cohomology class $b\alpha+c\beta \in H^2(M_1;\mathbb{Z})$, we
 can consider the Whitney sum of a complex line bundle over $M_1$ whose Chern class is $b\alpha+c\beta \in H^2(M_1;\mathbb{Z})$ and a trivial
 complex line bundle over $M_1$ and its projectivization, whose fiber $F_2$ is diffeomorphic to $S^2$. By this step, we obtain
 various manifolds $M_2$. More precisely, the Pontrjagin class of
 the manifold $M_2$ is represented as $c(2b-ac)\alpha \beta \in H^4(M_2;\mathbb{Z})$ where
 the two classes $\alpha$ and $\beta$ are classes defined as the
 pull-backs of the two classes $\alpha$ and $\beta$ before by the projection of the bundle $M_2$ over $M_1$, respectively. If $c$ is even, then the restriction of
 the bundle $M_2$ over $M_1$ to the fiber $F_1$ of the bundle $M_1$ before is trivial and
 as a result, by applying Theorem \ref{thm:4}, we obtain a round fold map from $M_2$ into ${\mathbb{R}}^2$. As the resulting manifolds $M_2$, we can
 obtain various manifolds by the observation of the Pontrjagin classes before. 
\item
\label{ex:3.3}
 Let $M$ be the $3$-dimensional complex projective space ${\mathbb{C}P}^3$. $M$ is regarded as
the total space of an $S^2$-bundle over the $4$-dimensional standard sphere.
By this bundle structure together with Theorem \ref{thm:1} (\ref{thm:1.1}), we obtain a round fold
map from $M$ into $\mathbb{R}^4$. We consider complex vector bundles of dimension $k \geq 1$
over $M$ and their subbundles whose fibers are the unit spheres of dimension
$2k-1$. If the 1st Chern class of such a complex vector bundle is represented as
$2c$ for a class $c \in H^2(M;Z)$, then by applying Theorem \ref{thm:4}, we obtain a round
fold map from the total space of the mentioned subbundle into ${\mathbb{R}}^4$; set $k=2$
and $(F_1,F_2)=(S^2,S^{2k-1})$. As in the previous example, in the case where $k$ is
larger than $1$, we can consider the projectivizations of complex vector bundles
and obtain bundles whose fibers are diffeomorphic to the ($k-1$)-dimensional
complex projective space ${\mathbb{C}P}^{k-1}$. If the 1st Chern class of such a complex
vector bundle is represented as $2c$ for a class $c \in H^2(M;Z)$, then by applying
Theorem \ref{thm:4}, we obtain a round fold map from the total space of the mentioned
${\mathbb{C}P}^{k-1}$-bundle into ${\mathbb{R}}^4$; set $k=2$ and $(F_1, F_2) = (S^2,{\mathbb{C}P}^{k-1})$.
For precise studies on the isomorphism classes of such bundles and the topologies and the differentiable stuructures of the total spaces of such bundles, see \cite{choimasudasuh} and \cite{kurokisuh} for
example.
\end{enumerate}
\end{Ex}


\section{The topologies and the differentiable structures of manifolds admitting round fold maps having good algebraic and differential topological conditions}
\label{sec:4}
  
  In the previous section, we introduce round fold maps on homotopy spheres and construct round fold maps on manifolds regarded as the total spaces of smooth bundles
 over standard spheres in Theorems \ref{thm:1}-\ref{thm:4}. Furthermore, we characterize manifolds regarded as total spaces of smooth bundles
 over standard spheres in Theorems \ref{thm:1} and \ref{thm:3} by certain families of round fold maps.

 In this section, we construct a new round fold map on a manifold represented as a connected sum of two manifolds admitting round fold maps under appropriate
 conditions and conversely, we decompose a round fold map into two round fold maps so that the original source manifold is represented as a
 connected sum of the resulting two source manifolds. Last, by using these two operations, we characterize
 some families of manifolds by some families of round fold maps. In this section, we mainly consider round fold maps from manifolds whose dimensions
 are $m$ into ${\mathbb{R}}^n$ satisfying the relation $m \geq 2n$.

 First, we review a result on homology and homotopy groups of manifolds admitting
 round fold maps whose inverse images of regular values are always disjoint unions of homotopy spheres. Maps satisfying the assumption of the proposition often appear
 as important examples in discussions of this section.

\begin{Prop}[\cite{kitazawa}, \cite{kitazawa2}, \cite{kitazawa3} etc.]
\label{prop:1}
Let $M$ be a closed and connected manifold of dimension $m$. Let $f:M \rightarrow {\mathbb{R}}^n$ be a round fold
 map satisfying the following two and the relation $m>n \geq 2$ holds. 
\begin{enumerate}
\item For each regular value p, $f^{-1}(p)$ is a disjoint union of almost-spheres.
\item Indices of fold points are $0$ or $1$ .
\end{enumerate}
Then, the quotient map $q_f:M \rightarrow W_f$ induces an isomorphism
 of homotopy groups ${\pi}_k(M) \cong {\pi}_k(W_f)$ for $0 \leq k \leq m-n-1$.

 Furthermore, let $M$ be simply-connected and let the relation $m \geq 2n$ also hold. Let the inverse image of a point
 in a proper core of $f$ consist of $l \in \mathbb{N}$ connected components. Then, we have ${\pi}_k(M) \cong \{0\}$ for $0 \leq k \leq n-1$. We also have
 ${\pi}_n(M) \cong H_n(M;\mathbb{Z}) \cong {\mathbb{Z}}^{l-1}$ in the case where $m>2n$ holds
 and ${\pi}_n(M) \cong H_n(M;\mathbb{Z}) \cong {\mathbb{Z}}^{2(l-1)}$ in the case where $m=2n$ holds. \\
\end{Prop} 

\begin{Rem}
As an important class of fold maps, {\it simple} fold maps are defined as fold maps such that the maps obtained as the restriction maps of the quotient maps onto the Reebs to the singular sets are injective. Special generic maps, round fold maps and fold maps such that the maps obtained by the restrictions to the singular sets are embeddings, are simple. For the definition and
 several algebraic and differential topological properties, see \cite{saeki} and \cite{sakuma2} for example.

The former part of Proposition \ref{prop:2} is not only for round fold maps satisfying the given two restrictions
 but also for simple fold maps from $m$-dimensional closed and connected manifolds into $n$-dimensional manifolds without boundaries satisfying
 the restrictions under the condition that $m-n \geq 2$ is assumed or that $m-n=1$ and the orientability of the manifold $M$ are assumed. See the referred articles and also \cite{saekisuzuoka}.

As an advanced study, in \cite{kitazawa4}, for appropriate cases of functions, we obtain special generic maps into the plane or higher dimensional Euclidean spaces such that by composing the canonical projections to $\mathbb{R}$ to the special generic maps, the obtained functions are the given functions. 
\end{Rem} 
\subsection{Round fold maps on connected sums of two connected manifolds admitting round fold maps}
\label{subsec:4.1}
In this subsection, we construct a round fold map on a manifold represented as a connected sum of
 two closed and connected manifolds admitting round fold maps under additional conditions.

 We introduce an operation of constructing a new round fold map from given two round fold maps.

 Let $m,n \in \mathbb{N}$. Let $M_1$ and $M_2$ be
 closed and connected manifolds whose dimensions are $m$. 
Let there exist round fold maps $f_1:M_1 \rightarrow {\mathbb{R}}^n$
 and $f_2:M_2 \rightarrow {\mathbb{R}}^n$ such that the two following hold.
\begin{enumerate}
\item The inverse image of a point in a proper core $P_1$ of $f_1$ has a connected component diffeomorphic to $S^{m-n}$.
\item For the boundary $C$ of the unbounded connected component of ${\mathbb{R}}^n-{\rm Int} f_2(M_2)$ ($n \neq 1$) or the disjoint union of the two unbounded connected components of ${\mathbb{R}}^n-{\rm Int} f_2(M_2)$ ($n=1$) and a small
 closed tubular neighborhood $P_2$, ${f_2}^{-1}(P_2)$ is the total space of a trivial linear $D^{m-n+1}$-bundle over the connected
 component $\partial P_2 \bigcap f_2(M_2)$ of $\partial P_2$ such
 that the map $f_2 {\mid}_{{f_2}^{-1}(\partial P_2 \bigcap f_2(M))}$ gives a subbundle of the
 bundle ${f_2}^{-1}(P_2)$.
\end{enumerate}

Let $V_1$ be a connected component of ${f_1}^{-1}(P_1)$ such that $f {\mid}_{V_1}:V_1 \rightarrow P_1$ makes $V_1$ a trivial smooth $S^{m-n}$-bundle over
 $D^n$, which exists by the assumption on $f_1$, and $V_2:={f_2}^{-1}(P_2)$. $V_2$ is a small closed tubular neighborhood
 of ${f_2}^{-1}(C) \subset M_2$. $V_2$ is regarded as the total space of a trivial
 linear $D^{m-n+1}$-bundle from the latter assumption. Moreover, from the same assumption, $f_2 {\mid}_{\partial V_2}$ gives
 a subbundle of the bundle $V_2$.

 We can consider two maps $f_1 {\mid}_{M_1-{\rm Int} V_1}$ and $f_2 {\mid}_{M_2-{\rm Int} V_2}$ and glue
 them together by a diffeomorphism regarded as a bundle isomorphism between the bundles $\partial V_2$ and $\partial V_1$ inducing a diffeomorphism between the
 base spaces to obtain a new round fold map $f$ from a new manifold into ${\mathbb{R}}^n$.

\begin{figure}
\includegraphics[width=50mm]{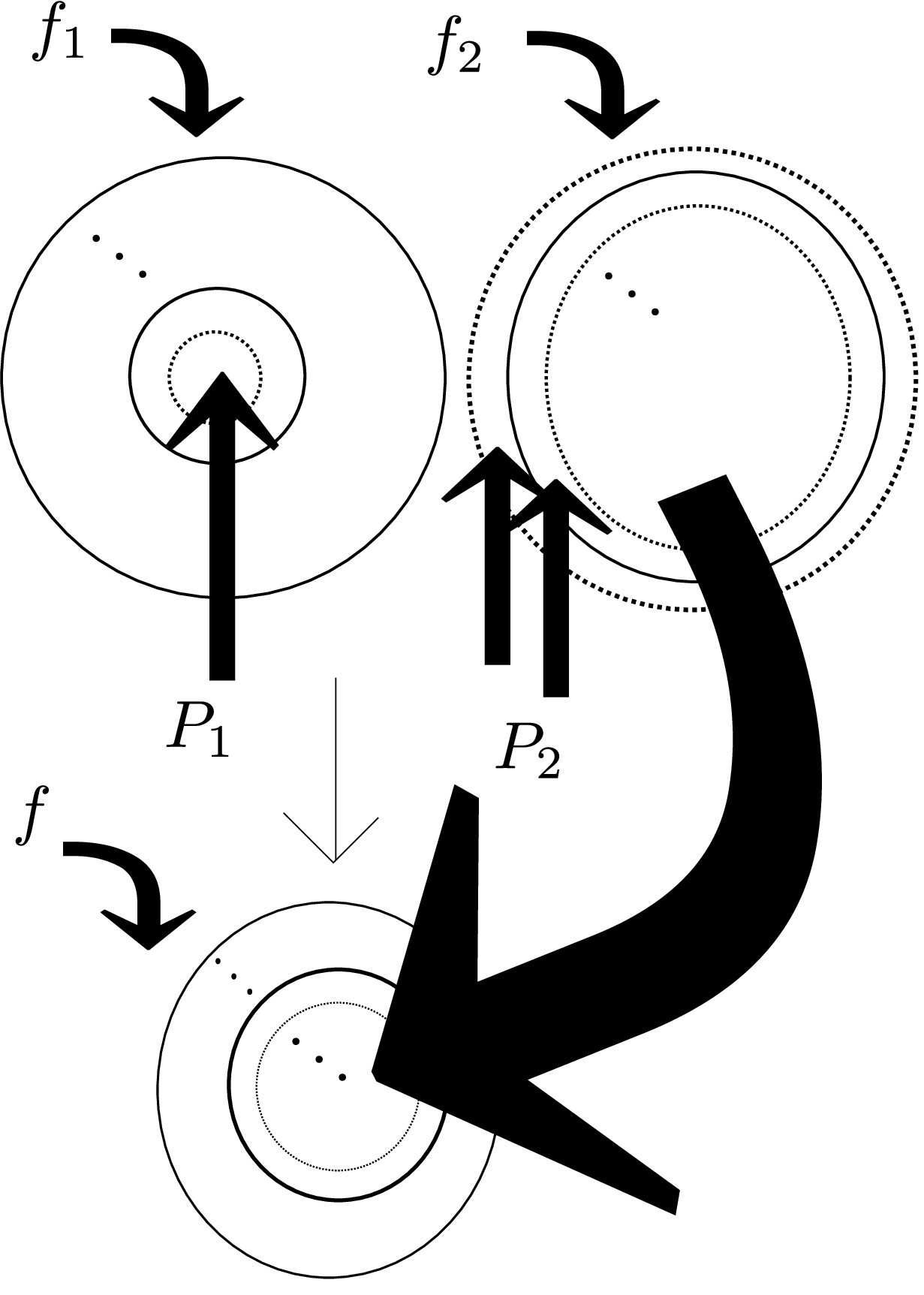}
\caption{A canonical combining operation to the pair $(f_1,f_2)$ to obtain $f$.}
\label{fig:4}
\end{figure}

 We call the operation a {\it canonical combining operation} to the pair $(f_1,f_2)$. See also FIGURE \ref{fig:4}.
 

\begin{Prop}
\label{prop:2}
Let $M_1$ and $M_2$ be closed and connected manifolds whose dimensions are $m$. 
Let there exist a round fold map into ${\mathbb{R}}^n$ {\rm (}$n \geq 1${\rm )} $f_1:M_1 \rightarrow {\mathbb{R}}^n$ such that the
 inverse image of a point in a proper core of $f_1$ has a connected component diffeomorphic to $S^{m-n}$
 and a round fold map $f_2:M_2 \rightarrow {\mathbb{R}}^n$ such that for the boundary $C$ of the unbounded connected component of ${\mathbb{R}}^n-{\rm Int} f_2(M_2)$ or the boundary $C$ of the disjoint union of the two unbounded connected component of ${\mathbb{R}}^n-{\rm Int} f_2(M_2)$, the inclusion
 of ${f_2}^{-1}(C)$ into $M_2$ is null-homotopic. We also assume that the relation $m \geq 2n$ holds.

 Then, on any manifold represented as a connected sum $M$ of $M_1$ and $M_2$, by a canonical combining operation to the pair $(f_1,f_2)$ we obtain
 a round fold map $f:M \rightarrow {\mathbb{R}}^n$.
\end{Prop}

\begin{proof}
Let $P_1$ be a proper core of $f_1$ and $P_2$ be a small closed tubular neighborhood of the connected component $C$ of $f_2(S(f_2))$. Let
 $V_1$ be a connected component of ${f_1}^{-1}(P_1)$ such that $f {\mid}_{V_1}:V_1 \rightarrow P_1$ gives a trivial smooth $S^{m-n}$-bundle over
 $D^n$, which exists by the assumption on $f_1$, and $V_2:={f_2}^{-1}(P_2)$. $V_2$ is a closed tubular neighborhood of ${f_2}^{-1}(C) \subset M_2$. Since $m \geq 2n=2(n-1)+2$ is assumed
 and the inclusion of ${f_2}^{-1}(C)$ into $M_2$ is assumed
 to be null-homotopic, $V_2$ is regarded as the total space of a trivial linear $D^{m-n+1}$-bundle. More precisely, $f_2 {\mid}_{\partial V_2}$ gives a subbundle of the bundle.

 Since $m \geq 2n=2(n-1)+2$ is assumed, we may regard that the following holds for each diffeomorphism
 $\Psi:\partial D^m \rightarrow \partial D^m$ extending to some diffeomorphism on $D^m$ or from $M_2-(M_2-D^m)$
 onto $M_1-(M_1-D^m)$ and for some diffeomorphism $\Phi:\partial V_2 \rightarrow \partial V_1$ regarded as a bundle
 isomorphism between the two trivial smooth $S^{m-n+1}$-bundles over $S^{n-1}$ inducing a diffeomorphism between the base
 spaces, where for two manifolds $X_1$ and $X_2$, $X_1 \cong X_2$ means
 that $X_1$ and $X_2$ are diffeomorphic. 

\begin{eqnarray*}
& & (M_1-{\rm Int} V_1) {\bigcup}_{\Phi} (M_2-{\rm Int} V_2) \\
& \cong & (M_1-{\rm Int} V_1) {\bigcup}_{\Phi} ((D^m-{\rm Int} V_2) \bigcup (M_2-{\rm Int} D^m)) \\
& \cong & (M_1-{\rm Int} V_1) {\bigcup}_{\Phi} ((S^m-({\rm Int} V_2 \sqcup {\rm Int} D^m)) {\bigcup}_{\Psi} (M_2-{\rm Int} D^m)) \\
& \cong & (M_1-{\rm Int} D^m) {\bigcup}_{\Psi} (M_2-{\rm Int} D^m)
\end{eqnarray*}

 This means that the resulting manifold $M$ is represented
 as a connected sum of the two manifolds $M_1$ and $M_2$ and that $M$ admits a round fold map $f:M \rightarrow {\mathbb{R}}^n$. More precisely, $f$
 is obtained by a canonical combining operation to the pair $(f_1,f_2)$. 
\end{proof}

\begin{Ex}
\label{ex:4}
Let $M_1$ and $M_2$ be closed and connected manifolds whose dimensions are $m$. Let there exist a round fold map into ${\mathbb{R}}^n$ {\rm (}$n \geq 1${\rm )} $f_1:M_1 \rightarrow {\mathbb{R}}^n$ such that the
 inverse image of a point in a proper core of $f_1$ has a connected component diffeomorphic to $S^{m-n}$. We also assume that the relation $m \geq 2n$ holds. If ${\pi}_{n-1}(M_2) \cong \{0\}$ holds and a round fold map $f_2:M \rightarrow {\mathbb{R}}^n$ exists, then
 the pair of the maps $f_1$ and $f_2$ satisfies the assumption of
 Proposition \ref{prop:2}. For example, if $M_2$ is simply connected and admits a round fold map $f_2:M_2 \rightarrow {\mathbb{R}}^n$ satisfying the assumption of Proposition \ref{prop:1}, then the pair of the two maps satisfies the assumption of
 Proposition \ref{prop:2}.   
\end{Ex}

\begin{Cor}
\label{cor:1}
Let $M_1$ be a closed and connected manifold of dimension $m$, let $k>0$ be an integer and let $\{M_{2,j}\}_{j=1}^k$ be a family of closed and connected manifolds of dimension $m$ whose {\rm (}$n-1${\rm )}-th homotopy groups vanish. 
Let there exist a round fold map $f_1:M_1 \rightarrow {\mathbb{R}}^n$ such that the
 inverse image of a point in a proper core of $f_1$ has just $k>0$ connected components diffeomorphic to $S^{m-n}$
 and a family of $k$ round fold maps $\{f_{2,k}:M_{2,k} \rightarrow {\mathbb{R}}^n\}$. We also assume that the relation $m \geq 2n$ holds.

 Then, on any manifold represented as a connected sum $M$ of $k+1$ manifolds $M_1$, $M_{2,1}$ $\cdots$ and $M_{2,k}$, by a finite iteration of canonical combining operations, starting from the pair $(f_1,f_{2,1})$, we obtain
 a round fold map $f:M \rightarrow {\mathbb{R}}^n$.
\end{Cor}

By applying this example or this corollary, we have the following theorem in the case
 where the map is from a $5$-dimensional manifold into the plane.

\begin{Thm}
\label{thm:5}
Let $M$ be a closed and simply-connected manifold of dimension $5$. Then, $M$ admits a
 round fold map $f:M \rightarrow {\mathbb{R}}^2$ satisfying the assumption of Proposition \ref{prop:1} if and only if $M$ is represented
 as a connected sum of a finite number of manifolds regarded as the total spaces of smooth $S^3$-bundles over $S^2$.
\end{Thm}

For the proof of this theorem, we need the following proposition.

\begin{Prop}[Barden, \cite{barden}]
\label{prop:3}
Let $M$ be a closed and simply-connected manifold of dimension $5$. Then, $H_2(M;\mathbb{Z})$ is torsion-free if and only if $M$ is
 represented as a connected sum of a finite number of manifolds regarded as the total spaces
 of smooth $S^3$-bundles over $S^2$.
\end{Prop}

\begin{proof}[Proof of Theorem \ref{thm:5}]
By Theorem \ref{thm:1} (\ref{thm:1.1}) or Example \ref{ex:2} and by Proposition \ref{prop:2} together with Example \ref{ex:4}, a manifold represented as a connected sum of a finite number of manifolds regarded as the total spaces of smooth $S^3$-bundles
 over $S^2$ admits a round fold map into ${\mathbb{R}}^2$ satisfying the assumption of Proposition \ref{prop:1}. This completes the proof of the "if" part. 

 Conversely, we explain about the "only if" part. We can know this part simply from the result of Proposition \ref{prop:1}. However, we show this by seeing the topology of the Reeb space to present an important polyhedron used later. For the argument here, see also the proof of Proposition \ref{prop:1} (Theorem 3 of \cite{kitazawa3}), section 3 of \cite{kitazawa} etc..
 We define the following objects.
\begin{enumerate}
\item $A$ is a disjoint union of finite copies of $D^n$ and $B:=S^{n-1} \times L$, where $L$ is a compact and connected graph with no loops (FIGURE \ref{fig:5}).
\item $\psi:S^{n-1} \times \Lambda \rightarrow \partial A$ is a ${\rm PL}$ homeomorphism, where $\Lambda$ is a set consisting of a finite number of degree $1$ vertices of the graph $L$.
\end{enumerate}
\begin{figure}
\includegraphics[width=50mm]{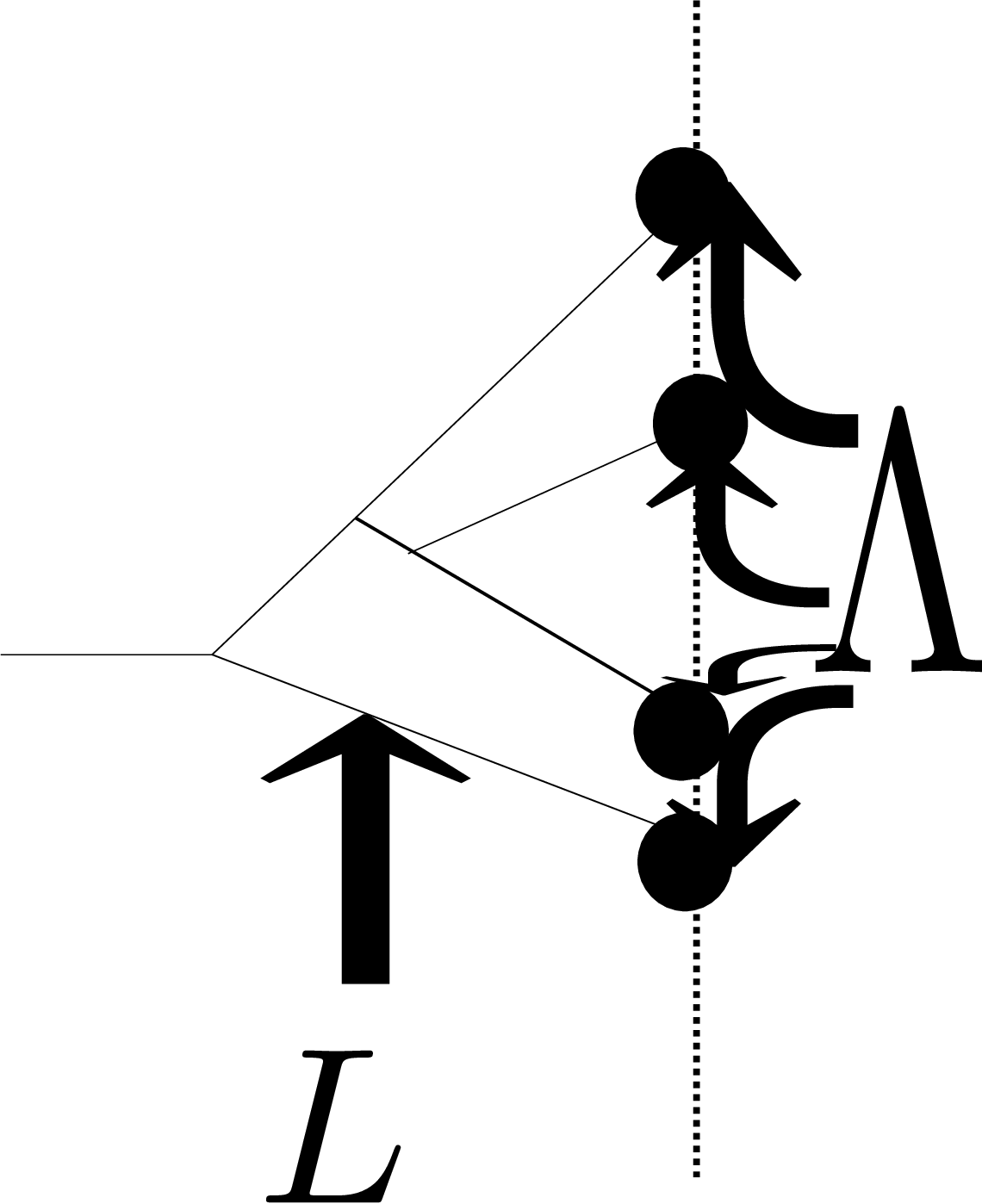}
\caption{The graph $L$.}
\label{fig:5}
\end{figure}
 
Then, for a round fold map $f$ from a closed and simply-connected manifold $M$ of dimension $5$ into ${\mathbb{R}}^2$ satisfying the
 assumption of Proposition \ref{prop:1}, the Reeb space $W_f$ is simple homotopy equivalent to a polyhedron given by attaching $B$ to $A$ by a homeomorphism $\phi$.
 It is regarded as a bouquet of finite copies of $S^2$ in the viewpoint of the (simple) homotopy theory. 

 From Proposition \ref{prop:1}, ${\pi}_1(M) \cong \{0\}$ holds and ${\pi}_2(M) \cong H_2(M;\mathbb{Z})$ is torsion-free. \\

 Hence, a closed and simply-connected manifold $M$ of dimension $5$ such that $H_2(M;\mathbb{Z})$ is not torsion-free does not admit a round fold map into ${\mathbb{R}}^2$ satisfying the assumption of Proposition \ref{prop:1}. From
 Proposition \ref{prop:3}, $H_2(M;\mathbb{Z})$ is not torsion-free if $M$ is not represented as a connected sum of a finite number
 of manifolds regarded as the total spaces of smooth $S^3$-bundles over $S^2$. This
 completes the proof of both part of the theorem. 
\end{proof}
\begin{Rem}
 Theorem \ref{thm:5} states that $5$-dimensional closed simply-connected manifolds
 with torsion-free homology groups admit round fold maps into the plane, however, we do not know examples of $5$-dimensional
 closed simply-connected manifolds
 which are also rational homology spheres not being integral homology spheres, discussed in \cite{barden}, admitting round fold maps into the plane. All $5$-dimensional
closed manifolds admit stable fold maps into the plane by virtue of well-known fundamental theory: a closed manifold admits a stable fold map into the plane if and only if the Euler number of the manifold is even (\cite{thom}, \cite{whitney} etc.).
\end{Rem}
\subsection{Decompositions of a manifold admitting a round fold map into a connected sum of two manifolds}
\label{subsec:4.2}

Conversely, in this subsection, we decompose a manifold admitting a round fold map into a connected sum of two manifolds
 admitting round fold maps under appropriate conditions. 

 Let $M$ be a closed and connected manifold of dimension $m$. Let $f:M \rightarrow {\mathbb{R}}^n$ be a round fold map.

Let $C$ be a submanifold
 of ${\mathbb{R}}^n-f(S(f))$. Suppose that $C$ is diffeomorphic to $S^{n-1}$ and that $C$ satisfies 
 one of the following.
\begin{enumerate}
\item $n \neq 1$ holds and $C$ is a deformation retract of
 the closure of a connected component of ${\mathbb{R}}^n-f(S(f))$ diffeomorphic to $S^{n-1} \times (0,1)$.
\item $n=1$ holds and $C$ is a set consisting of just $2$ points in mutually distinct connected components of ${\mathbb{R}}^n-f(S(f))$ such that in a normal form of the function, the corresponding image of the disjoint union of the two connected components is represented as $(-l-1,-l) \sqcup (l,l+1)$ for a positive integer $l$.
\item  $C$ is in a proper core of $f$. 
\end{enumerate}
We
 denote the closure of the bounded domain of ${\mathbb{R}}^n-C$ by $R$ and
 for a connected component $\bar{M}$ of $f^{-1}(R)$, $f {\mid}_{\partial \bar{M}}:\partial \bar{M} \rightarrow C$
 makes $\partial \bar{M}$ a trivial smooth bundle over $C$ and the inverse image is a standard sphere.

We can glue the map $f {\mid}_{M-{\rm Int} \bar{M}}:M-{\rm Int} \bar{M} \rightarrow {\mathbb{R}}^n$ and the natural
 projection $p:R \times S^{m-n} \rightarrow R$ by a diffeomorphism ${\Phi}_1:\partial M \rightarrow \partial R \times S^{m-n}$ regarded as a bundle isomorphism between the trivial $S^{m-n}$-bundles and a
 diffeomorphism between the base spaces induced from the isomorphism to obtain a round fold map $f:M_1 \rightarrow {\mathbb{R}}^n$. 

We can glue $\bar{M}$ and $V:=S^{n-1} \times D^{m-n+1}$ on the boundaries by any bundle isomorphism ${\Phi}_2$ between the two trivial smooth bundles $V:=S^{n-1} \times \partial D^{m-n+1}$ and $\partial \bar{M}$ over the {\rm (}$n-1${\rm )}-dimensional
 standard spheres inducing
 a diffeomorphism between the base spaces to obtain a manifold $M_2$. 

Let $\tilde{f_0}:D^{m-n+1} \rightarrow [0,+\infty)$ be a Morse function satisfying the following three as introduced just after Definition \ref{def:3}.
\begin{enumerate}
\item $\tilde{f_0}$ is constant and $0$ on the boundary.
\item $\tilde{f_0}({\rm Int} D^{m-n+1}) \subset (0,+\infty)$ holds.
\item $\tilde{f_0}$ has just one singular point and it is maximal at this point.
\end{enumerate}

We can also glue the maps $f {\mid}_{\bar{M}}:\bar{M} \rightarrow R$ and the product of the Morse function $\tilde{f}$ and
 the identity map ${\rm id}_{S^{n-1}}$ on $S^{n-1}$ by the diffeomorphism ${\Phi}_2$ and
 a diffeomorphism between the base spaces induced from the isomorphism to obtain a round fold map $f_2:M_2 \rightarrow {\mathbb{R}}^n$. 

We call this operation of obtaining the pair of two round fold maps $f_1$ and $f_2$ a {\it canonical decomposing operation} to the map $f$ over $C$.


Canonical decomposing operations are regarded as extensions of several specific {\it R-operations} introduced in \cite{kobayashisaeki}, which are regarded as surgery
 operations on stable maps from closed manifolds whose dimensions are larger than $2$ into the plane; in the paper, the operations are introduced as
 surgery operations to {\it pseudo quotient maps}, which are regarded as continuous maps from closed manifolds whose dimensions are larger than $2$
 into $2$-dimensional polyhedra locally topologically regarded as maps obtained from stable maps
 from the closed manifolds into the plane by considering natural quotient maps into the Reeb spaces. 

Note also that in \cite{kobayashisaeki}, canonical combining operations, which are regarded as inverse operations of
 the operations, are not defined. 

We introduce and show Proposition \ref{prop:4}. Compare this with arguments on R-operations shown in \cite{kobayashisaeki}.

\begin{Prop}
\label{prop:4}
Let $M$ be a closed and connected manifold of dimension $m$. Let $f:M \rightarrow {\mathbb{R}}^n$ be a round
 fold map. Suppose also that the relation $m \geq 2n$ holds.

Let $C$ be a submanifold
 of ${\mathbb{R}}^n-f(S(f))$. Suppose that $C$ is diffeomorphic to $S^{n-1}$ and
that $C$ satisfies 
 one of the following.
\begin{enumerate}
\item $n \geq1$ holds and $C$ is a deformation retract of
 the closure of a connected component of ${\mathbb{R}}^n-f(S(f))$ diffeomorphic to $S^{n-1} \times (0,1)$.
\item $n=1$ holds and $C$ is a set consisting of just $2$ points in mutually distinct connected components of ${\mathbb{R}}^n-f(S(f))$ such that in a normal form of the function, the corresponding image of the disjoint union of the two connected components is represented as $(-l-1,-l) \sqcup (l,l+1)$ for a positive integer $l$.
\item  $C$ is in a proper core of $f$. 
\end{enumerate}
We denote the closure of the bounded domain of ${\mathbb{R}}^n-C$ by $R$ and
 let there exist a connected component $\bar{M}$ of $f^{-1}(R)$ satisfying the following two.  
\begin{enumerate}
\item $f {\mid}_{\partial \bar{M}}:\partial \bar{M} \rightarrow C$ gives a trivial smooth bundle over $C$ and the fiber is a standard sphere.
\item Let $p$ be a point in $D^{m-n+1}$. If we glue $\bar{M}$ and $V:=S^{n-1} \times D^{m-n+1}$ on the boundaries by a bundle isomorphism $\Phi$ between
 the two trivial smooth bundles $V:=S^{n-1} \times \partial D^{m-n+1}$ and $\partial \bar{M}$ over the {\rm (}$n-1${\rm )}-dimensional standard spheres inducing
 a diffeomorphism between the base spaces, then the natural inclusion
 $S^{n-1} \times \{p\} (\subset V=S^{n-1} \times D^{m-n+1}) \subset \bar{M} {\bigcup}_{\Phi} V$ is null-homotopic.
\end{enumerate}

 Then, $M$ is represented
 as a connected sum of two connected manifolds $M_1$ and $M_2:=\bar{M} {\bigcup}_{\Phi} V$ such that $M_i$ admits a
 round fold map $f_i:M_i \rightarrow {\mathbb{R}}^n$ {\rm (}$i=1,2${\rm )} and that the pair $(f_1,f_2)$ is obtained by a canonical decomposing operation to the map $f$ over $C$.
\end{Prop}

\begin{proof}
Let $M_1$ be a closed manifold
 given by gluing $M-{\rm Int} \bar{M}$ and $V^{\prime}:=D^n \times S^{m-n}$ by some diffeomorphism
 ${\Phi}^{\prime}:\partial D^n \times S^{m-n} \rightarrow \partial \bar{M}$
 regarded as a bundle isomorphism between the natural two trivial smooth $S^{m-n}$-bundles inducing a diffeomorphism between the base spaces.

Since $m \geq 2n=2(n-1)+2$ is assumed and the natural inclusion
 $S^{n-1} \times \{p\} \subset V=S^{n-1} \times D^{m-n+1} \subset \bar{M} {\bigcup}_{\Phi} V$ is assumed to be null-homotopic, we may regard that the following holds for an orientation reversing diffeomorphism
 $\Psi:\partial D^m \rightarrow \partial D^m$ extending to a diffeomorphism on $D^m$ or from $M_2-(M_2-D^m)$ onto $M_1-(M_1-D^m)$, where for two manifolds $X_1$ and $X_2$, $X_1 \cong X_2$ means
 that $X_1$ and $X_2$ are diffeomorphic. 

\begin{eqnarray*}
& & (M-{\rm Int} \bar{M}) {\bigcup} \bar{M} \\
& \cong & (M_1-{\rm Int} {V}^{\prime}) {\bigcup} (M_2-{\rm Int} V) \\
& \cong & (M_1-{\rm Int} {V}^{\prime}) {\bigcup} ((D^m-{\rm Int} V) \bigcup (M_2-{\rm Int} D^m)) \\
& \cong & (M_1-{\rm Int} {V}^{\prime}) {\bigcup} ((S^m-({\rm Int} V \sqcup {\rm Int} D^m)) {\bigcup}_{\Psi} (M_2-{\rm Int} D^m)) \\
& \cong & (M_1-{\rm Int} D^m) {\bigcup}_{\Psi} (M_2-{\rm Int} D^m)
\end{eqnarray*}     

This means that $M$
 is represented as a connected sum of the two connected manifolds $M_1$ and $M_2$ and that $M_i$ admits
 a round fold map $f_i:M_i \rightarrow {\mathbb{R}}^n$ satisfying the mentioned
 condition ($i=1,2$). 
\end{proof}

\subsection{Explicit cases}
\label{subsec:4.3}
In this subsection, by applying the previous results and their proofs, we study the topologies and the differentiable structures of manifolds admitting round fold maps under appropriate differential topological constraints. More precisely, as a work, we characterize the
 family of all $m$-dimensional manifolds represented as connected sums of finite numbers
 of manifolds regarded as the total spaces of smooth $S^{m-n}$ bundles over $S^n$ by a class of round fold maps
 including a map in Example \ref{ex:2} under the assumption that the relation $m \geq 2n$ hold. 

\begin{Thm}
\label{thm:6}
Let $M$ be a closed and connected manifold of dimension $m$. Suppose that there exists
 a round fold map $f:M \rightarrow {\mathbb{R}}^n$. Suppose aso that the relation $m \geq 2n$ holds.

Let $N((S(f)))$ be a small closed tubular neighborhood
 of $f(S(f))$ and for any connected component $R$
 of ${\mathbb{R}}^n-{\rm Int} N(f(S(f)))$, $f {\mid}_{{f}^{-1}(R)}:f^{-1}(R) \rightarrow R$
 gives a trivial smooth bundle whose fiber is a disjoint union of standard spheres. Assume also that the number of connected
 components of the inverse image of a point in a proper core of $f$ equals the number of connected components of $S(f)$ and that the numbers are both larger than $2$.

Then, $M$ is represented as a connected sum of closed and connected manifolds each of which admits a
 round fold map. Furthermore, the following two hold.
\begin{enumerate}
\item The singular set of each round fold map above has just two connected components in the case where $n \neq 1$ is assumed and just four points in the case where $n=1$ is assumed and the inverse image of each point in its proper core
 consists of two standard spheres. 
\item Each round fold map just before is obtained by applying canonical decomposing operations to the round fold maps inductively starting from the map $f$.  
\end{enumerate}

Conversely, such a manifold
 admits a round fold map into ${\mathbb{R}}^n$ satisfying the assumption mentioned above. More precisely, it is obtained by
 applying canonical combining operations to the pairs of round fold maps inductively, starting
 with a family of round fold maps such that the singular set of each map consists of two connected components in the case $n \neq 1$ and just four points in the case $n=1$ and that the inverse image of each point in a
 proper core of each map consists of two standard spheres. 
\end{Thm}

\begin{proof}
All but one (the $n \geq 1$ case) or two (the $n=1$ case) connected components of the singular set $S(f)$ of $f$ consist of fold points whose indices are $1$ by the assumption that the number of connected
 components of the inverse image of a point in a proper core of $f$ equals the number of connected components of $S(f)$ and by the theory of attachments
 of handles. Furthermore, $f$ satisfies the assumption of Proposition \ref{prop:1}.

We define the following objects to represent the Reeb space $W_f$ of $f$., appearing also in the proof of Theorem \ref{thm:5}.
\begin{enumerate}
\item $A$ is a disjoint union of finite copies of $D^n$ and $B:=S^{n-1} \times L$, where $L$ is a compact and connected graph with no loops.
\item $\psi:S^{n-1} \times \Lambda \rightarrow \partial A$ is a ${\rm PL}$ homeomorphism, where $\Lambda$ is a set consisting
 of a finite number of degree $1$ vertices of the graph $L$.
\end{enumerate}
 Then, the Reeb space
 $W_f$ is ${\rm PL}$ homeomorphic to a polyhedron obtained by attaching $B$ to $A$ by a ${\rm PL}$ homeomorphism $\psi$. Note
 that ${\pi}_{n-1}(M) \cong {\pi}_{n-1}(W_f) \cong \{0\}$ holds from Proposition \ref{prop:1}.

By Proposition \ref{prop:4}, if the
 singular set of $f$ consists of more than two connected components in the case where $n \neq 1$ and more than four points in the case where $n=1$, then
 by a canonical decomposing operation, we obtain two round fold
 maps; one is a map satisfying the assumption of the theorem such that the number of connected components of the singular set of the map is smaller than that of the original map $f$ by one in the case where $n \neq 1$ and by two in the case where $n=1$ and the other map satisfies the assumption of this
 theorem satisfying that the number of connected components of the singular set of the map is two in the case where $n \neq 1$ and four in the case where $n=1$. Of course, $M$ is represented as a connected sum of the resulting source manifolds. It also follows that $M$ is represented as
 a connected sum
 of closed and connected manifolds admitting round fold maps whose singular sets consist of two or four connected components
 satisfying that the inverse images of points in proper cores consist of two standard spheres.

Conversely, if $M$ is represented as a connected sum of closed and connected manifolds
 admitting round fold maps whose singular sets always consist of two ($n \neq 1$) or four ($n=1$) connected components such that the inverse images of points in proper cores
 are always disjoint unions of two standard spheres, then, by virtue of Proposition \ref{prop:2} together with Proposition \ref{prop:1} or Example \ref{ex:3}, by
 using canonical combining operations inductively to the pairs of
 previous maps whose singular sets consist of
 two or four connected components, we can obtain a round fold map satisfying the assumption. 
\end{proof}

  Let ${\Theta}_{k}$ be the {\it h-cobordism group} of oriented homotopy spheres of dimension $k \geq 2$. It follows
 easily that the set
 of all classes of ${\Theta}_{k_1}$ consisting of oriented homotopy spheres admitting round fold maps with connected singular sets
 into ${\mathbb{R}}^{k_2}$ ($k_1 \geq k_2 \geq 2$) is a subgroup of ${\Theta}_{k_1}$ (note that the maps are special generic). In fact, by considering a {\it connected
 sum} of given two round fold maps with connected singular sets, we easily see that this set is a subgroup (for a {\it connected sum} of given two special generic maps into an Euclidean space,
 see section 5 of \cite{saeki2} for example). 

\begin{Def}
\label{def:4}
We denote the subgroup by ${\Theta}_{(k_1,k_2)} \subset {\Theta}_{k_1}$ and call it the {\it $(k_1,k_2)$ round special generic group}.  
\end{Def}

 We have the following two theorems.

\begin{Thm}
\label{thm:7}
Let $m,n>0$ be integers satisfying the relation $m \geq 2n$. Let $M$ be a closed and connected manifold
 of dimension $m$. Let $E$ be a manifold of dimension $m-n+1$ satisfying at least one of the following two.
\begin{enumerate}
\item The boundary $\partial E$ is non-empty and includes a {\rm (}$n-1${\rm )}-connected manifold.
\item The boundary $\partial E$ is non-empty and ${\pi}_{n-1}(E) \cong \{0\}$ holds.
\end{enumerate}

Let $M$ admit
 a round fold map $f:M \rightarrow {\mathbb{R}}^n$ such
 that the following two hold.
\begin{enumerate}
\item There exists a positive integer $l>1$ and the singular set $S(f)$ consists of $l$ connected components in the case $n \neq 1$ and $2l$ connected components in the case $n=1$.
\item Let $f_0:M \rightarrow {\mathbb{R}}^n$ be a normal form of the round fold map $f$. The composition of the map 
 ${f_0} {\mid}_{{f_0}^{-1}({D^n}_{l-\frac{1}{2}}-{\rm Int} {D^n}_{\frac{1}{2}})}:{f_0}^{-1}({D^n}_{l-\frac{1}{2}}-{\rm Int} {D^n}_{\frac{1}{2}}) \rightarrow {D^n}_{l-\frac{1}{2}}-{\rm Int} {D^n}_{\frac{1}{2}}$
 and the map ${\pi}_{P}(x):=\frac{1}{2}\frac{x}{|x|}$ {\rm (}$x \in {D^n}_{l-\frac{1}{2}}-{\rm Int} {D^n}_{\frac{1}{2}}${\rm )}
 defined just after Definition \ref{def:2} only in the case $n \neq 1$, we can define in the case $n=1$ of course, defines a bundle structure by Ehresmann's fibration theorem. The
 bundle is a trivial smooth bundle and the fiber is diffeomorphic to the manifold $E$ with the interior of a standard closed disc $D^{m-n+1}$ smoothly embedded
 in the interior ${\rm Int} E$ removed.
\end{enumerate}
 Then, by a canonical decomposing operation to $f$ over a sphere $C$ in the connected component of the regular value set of $f$ whose boundary contains the boundary $\partial f(M)$, we have a pair of
 round fold maps such that one of the maps is a round fold map from
 a manifold in a class of ${\Theta}_{(m,n)}$ into ${\mathbb{R}}^n$ whose singular set
 is connected, that the
 other is a $C^{\infty}$ trivial round fold map and that $M$ is
 represented as a connected sum of the resulting two source manifolds. In
 addition, we can construct the latter
 map so that the inverse image of its axis is diffeomorphic to $E$ and that the number of connected components of
 its singular set is $l$.

Especially, if $\partial E$ has a connected component diffeomorphic to the standard sphere $S^{m-n}$, then $M$ admits
 a $C^{\infty}$ trivial round fold map satisfying the following conditions in the case $n \neq 1$.

\begin{enumerate}
\item The singular set consists of $l$ connected components.
\item The inverse image of an axis of the map is diffeomorphic to $E$.
\end{enumerate}
Moreover, $M$ admits
 a round fold map with $2l$ singular points and the inverse image of an axis of the map is diffeomorphic to $E$.
\end{Thm}

\begin{proof}
 From the assumption on the topology of manifold $E$, the two assumptions on the map $f$
 and Proposition \ref{prop:4}, the first result follows. We obtain
 the second result by applying Proposition \ref{prop:2}. More precisely, we can
 use a canonical combining operation to the pair of the round fold map
 whose singular set consists of $l$ connected components and the round fold map
 whose singular set is connected obtained by the canonical decomposing
 operation mentioned in the statement
 to obtain a new round fold map from $M$. The resulting round fold map on $M$ is $C^{\infty}$ trivial
 by the method of the construction. 
\end{proof}

\begin{Thm}
\label{thm:8}
Let $m, n>0$ be integers satisfying the relation $m \geq 2n$. Let $M$ be a closed and connected manifold of dimension $m$ and
 let $\Sigma$ be an almost-sphere of dimension $m-n$. Then, the following two are equivalent.
\begin{enumerate}
\item
\label{thm:8.1}
 A round fold map $f:M \rightarrow {\mathbb{R}}^n$ satisfying the following two exists.
\begin{enumerate}
\item $S(f)$ consists of two connected components in the case $n \neq 1$ and four connected components in the case $n=1$.
\item Let $f_0:M \rightarrow {\mathbb{R}}^n$ be a normal form of the round fold map $f$. The composition of the map 
 ${f_0} {\mid}_{{f_0}^{-1}({D^n}_{\frac{3}{2}}-{\rm Int} {D^n}_{\frac{1}{2}})}:{f_0}^{-1}({D^n}_{\frac{3}{2}}-{\rm Int} {D^n}_{\frac{1}{2}}) \rightarrow {D^n}_{\frac{3}{2}}-{\rm Int} {D^n}_{\frac{1}{2}}$
 and the map ${\pi}_{P}(x):=\frac{1}{2}\frac{x}{|x|}$ {\rm (}$x \in {D^n}_{\frac{3}{2}}-{\rm Int} {D^n}_{\frac{1}{2}}${\rm )} defines a bundle structure by Ehresmann's fibration theorem. The
 bundle is a trivial smooth bundle and the fiber is diffeomorphic to the
 cylinder $\Sigma \times [-1,1]$ with the interior of a standard closed disc $D^{m-n+1}$ smoothly embedded in
 the interior of the manifold removed.
\end{enumerate} 
\item
\label{thm:8.2}
 $M$ is represented as a connected sum of a manifold in a class of ${\Theta}_{(m,n)}$ and a manifold regarded as the total space of a smooth bundle over $S^n$ whose
 fiber is diffeomorphic to the almost-sphere $\Sigma$. If $\Sigma$ is diffeomorphic to the standard sphere $S^{m-n}$, then $M$ is
 a manifold regarded as the total space of a smooth bundle over $S^n$ whose fiber is diffeomorphic to the standard sphere $S^{m-n}$. 
\end{enumerate} 
\end{Thm}

\begin{proof}
Assume that a round fold map $f:M \rightarrow {\mathbb{R}}^n$ satisfying the condition (\ref{thm:8.1}) exists. Note
 that ${\pi}_{n-1}(\Sigma \times [-1,1]) \cong {\pi}_{n-1}(\Sigma) \cong \{0\}$ holds since the condition $m \geq 2n$ is assumed.  By virtue of Theorem \ref{thm:7},
 we can represent $M$ as a connected sum of two closed and connected manifolds $M_1$ and $M_2$ such that the following two hold. 

\begin{enumerate}
\item $M_1$ admits a round fold map $f_1:M_1 \rightarrow {\mathbb{R}}^n$ such that $S(f_1)$ is connected in the case $n \neq 1$ and a set with just $2$ points in the case $n=1$.
\item $M_2$ admits a $C^{\infty}$ trivial round fold map $f_2:M_2 \rightarrow {\mathbb{R}}^n$ such that $S(f_2)$ consists
 of two connected components and that the inverse image of an axis of $f_2$ is diffeomorphic to $\Sigma \times [-1,1]$ in the case $n \neq 1$ and a round fold map $f_2:M_2 \rightarrow {\mathbb{R}}^n$ such that $S(f_2)$ consists
 of four points in the case $n=1$. $f_2$ satisfies the assumption of Proposition \ref{prop:1}. 
\end{enumerate}

 Furthermore, the pair $(f_1,f_2)$ is obtained by a canonical decomposing operation to the map $f$. Then, $M$ is represented as a connected sum of a manifold in a
 class of ${\Theta}_{(m,n)}$ and a manifold regarded as the total space of a smooth bundle over $S^n$ whose fiber
 is diffeomorphic to the almost-sphere $\Sigma$ by Theorem \ref{thm:1} (\ref{thm:1.2}). If $\Sigma$ is diffeomorphic to $S^{m-n}$, then, by virtue of the second statement of Theorem \ref{thm:7}, $M$ admits a $C^{\infty}$ trivial round
 fold map into ${\mathbb{R}}^n$ such that the singular set consists
 of two connected components in the case $n \neq 1$ and four points in the case $n=1$ and that the inverse image of an axis of the map is diffeomorphic to $S^{m-n} \times [-1,1]$ and $M$ is a manifold regarded as the total space of
 a smooth bundle over $S^n$ whose fiber is diffeomorphic to the standard sphere $S^{m-n}$ by Theorem \ref{thm:1} (\ref{thm:1.2}).

Conversely, if $M$ is such a manifold, then by Proposition \ref{prop:2} or Example \ref{ex:4}
 together with Theorem \ref{thm:1} (\ref{thm:1.1}) and Example \ref{ex:2} (\ref{ex:2.1}), $M$ admits a
 round fold map $f:M \rightarrow {\mathbb{R}}^n$ satisfying the condition (\ref{thm:8.1}).

This completes the proof.
\end{proof}

Now we have the following theorem.

\begin{Thm}
\label{thm:9}
Let $m,n>0$ be integers satisfying the relation $m \geq 2n$. Let $M$ be a closed and connected manifold of dimension $m$. Then
 the following two are equivalent.

\begin{enumerate}
\item
\label{thm:9.1}
 A round fold map $f:M \rightarrow {\mathbb{R}}^n$ satisfying the following two exist.
\begin{enumerate}
\item Inverse images of regular values of are disjoint unions of standard spheres and the number of connected components of the inverse images of a point in a proper core of $f$
 and the number of connected components of $S(f)$ are $l>1$ in the case $n \geq 1$ and $l>2$ in the case $n=1$. 
\item Let $f_0:M \rightarrow {\mathbb{R}}^n$ be a normal form of the round fold map $f$. The composition of the map 
 ${f_0} {\mid}_{{f_0}^{-1}({D^n}_{k+\frac{1}{2}}-{\rm Int} {D^n}_{k-\frac{1}{2}})}:{f_0}^{-1}({D^n}_{k+\frac{1}{2}}-{\rm Int} {D^n}_{k-\frac{1}{2}}) \rightarrow {D^n}_{k+\frac{1}{2}}-{\rm Int} {D^n}_{k-\frac{1}{2}}$
 and the map ${\pi}_{P}(x):=\frac{1}{2}\frac{x}{|x|}$ defines a bundle structure by Ehresmann's fibration theorem for
 any integer $1 \leq k \leq l-1$. The
 bundles are all trivial smooth bundles and the fibers are disjoint unions of finite copies of the cylinder $S^{m-n} \times [-1,1]$ with
 the interior of a standard closed disc $D^{m-n+1}$ smoothly embedded in
 the interior of the manifold removed.
\end{enumerate}
\item
\label{thm:9.2}
$M$ is represented as a connected sum of $l-1$ manifolds regarded as the total spaces of smooth $S^{m-n}$-bundles over $S^n$.
\end{enumerate}
\end{Thm}

\begin{figure}
\includegraphics[width=50mm]{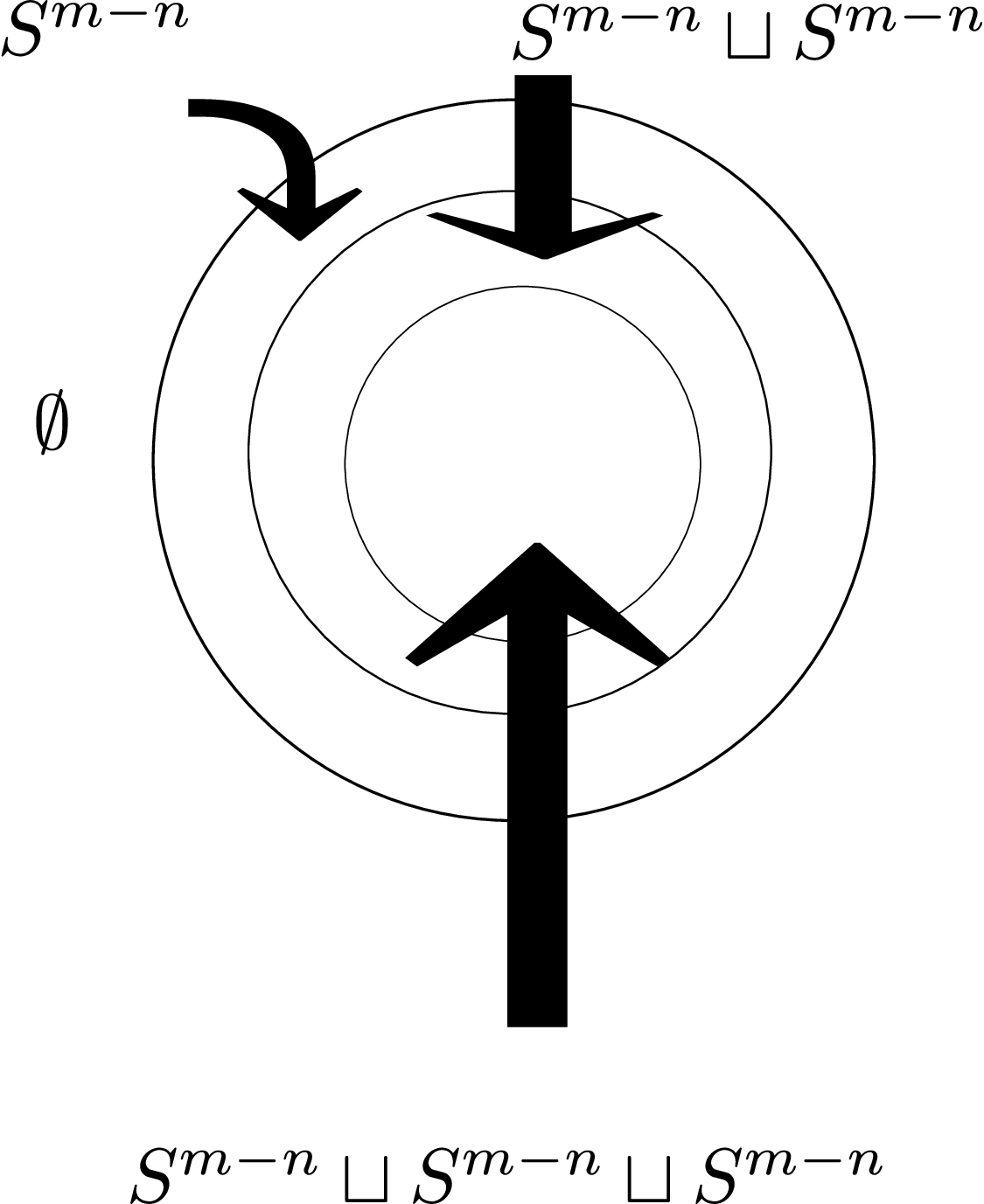}
\caption{The image (and the singular value set) of a round fold map into ${\mathbb{R}}^n$ in Theorem \ref{thm:9} : manifolds and the empty set represent inverse images of regular values.}
\label{fig:6}
\end{figure}

\begin{proof}
Assume that a round fold map $f:M \rightarrow {\mathbb{R}}^n$ satisfying the condition (\ref{thm:9.1}) exists.

By Theorems \ref{thm:6} and \ref{thm:8}, $M$ is represented as a
 connected sum of a finite number of closed manifolds
 admitting round fold maps whose singular sets consist of two connected components and each manifold is a manifold regarded as the total bundle of
 a smooth $S^{m-n}$-bundle over $S^n$.

Conversely, assume that the condition (\ref{thm:9.2}) holds. By Theorem \ref{thm:1} (\ref{thm:1.1}) and Example \ref{ex:2} (\ref{ex:2.1}), each manifold appearing as an ingredient of the connected sum admits a $C^{\infty}$ trivial round fold map whose singular set
 consists of two connected components such that the inverse image of an axis of the map is diffeomorphic to the cylinder $S^{m-n} \times [-1,1]$. By
 Proposition \ref{prop:2} or Example \ref{ex:4}, a manifold $M$ represented as a connected
 sum of a finite number of manifolds regarded as the total spaces of smooth $S^{m-n}$-bundles over
 $S^n$ admits a round fold map $f:M \rightarrow {\mathbb{R}}^n$ satisfying
 the condition (\ref{thm:9.1}).

This completes the proof.
\end{proof}
For maps used here, see also FIGURE \ref{fig:6}.
Last, we show several theorems.

\begin{Thm}
\label{thm:10}
Let $m, n>0$ be integers satisfying the relation $m \geq 2n$ and let $M$ be a closed and connected manifold of dimension $m$. Let $F$ be a closed and
 {\rm (}n-1{\rm )}-connected manifold of dimension $m-n$. Then 
 the following two are equivalent.

\begin{enumerate}
\item
\label{thm:10.1}
 A round fold map $f:M \rightarrow {\mathbb{R}}^n$ satisfying the following two exists.
\begin{enumerate}
\item The number of connected component $l$ of the singular set $S(f)$ is larger than $2$ in the case $n \neq 1$ and larger than $4$ in the case $n=1$.
\item Let $l^{\prime}$ be an integer satisfying $1 \leq l^{\prime}<l-1$. Let $f_0:M \rightarrow {\mathbb{R}}^n$ be a normal form of the round fold map $f$. The composition of the map 
 ${f_0} {\mid}_{{f_0}^{-1}({D^n}_{l^{\prime}+\frac{1}{2}}-{\rm Int} {D^n}_{\frac{1}{2}})}:{f_0}^{-1}({D^n}_{l^{\prime}+\frac{1}{2}}-{\rm Int} {D^n}_{\frac{1}{2}}) \rightarrow {D^n}_{l^{\prime}+\frac{1}{2}}-{\rm Int} {D^n}_{\frac{1}{2}}$
 and the map ${\pi}_{P}(x):=\frac{1}{2}\frac{x}{|x|}$ and the composition of the map 
 ${f_0} {\mid}_{{f_0}^{-1}({D^n}_{l-\frac{1}{2}}-{\rm Int} {D^n}_{l^{\prime}+\frac{1}{2}})}:{f_0}^{-1}({D^n}_{l-\frac{1}{2}}-{\rm Int} {D^n}_{l^{\prime}+\frac{1}{2}}) \rightarrow {D^n}_{l-\frac{1}{2}}-{\rm Int} {D^n}_{l^{\prime}+\frac{1}{2}}$
 and the map ${\pi}_{P}(x):=\frac{1}{2}\frac{x}{|x|}$ define bundle structures by Ehresmann's fibration theorem. Both bundles are
 trivial smooth bundles. The fibers are the cylinder $S^{m-n} \times [-1,1]$ with
 the interior of a standard closed disc $D^{m-n+1}$ smoothly embedded in
 the interior of the manifold removed and a manifold represented as the disjoint union
 of a manifold diffeomorphic to $F \times [-1,1]$ with the interior of an {\rm (}$m-n+1${\rm )}-dimensional standard
 closed disc smoothly embedded in the interior removed
 and the cylinder $S^{m-n} \times [-1,1]$, respectively. 
\end{enumerate}
\item
\label{thm:10.2}
$M$ is represented as a connected sum of a manifold regarded as the total spaces of smooth $S^{m-n}$-bundles over $S^n$ and a manifold regarded as the total space of a smooth $F$-bundle over $S^n$.
\end{enumerate}
\end{Thm}

\begin{proof}
Assume that a round fold map $f:M \rightarrow {\mathbb{R}}^n$ satisfying the condition (\ref{thm:10.1}) exists. These assumptions
 enable us to apply Proposition \ref{prop:4} and by a canonical decomposing operation, we obtain a
 pair of round fold maps. More precisely, one of the round fold maps is a $C^{\infty}$-trivial round fold map as in Theorem \ref{thm:1} such that the inverse image of a point in a proper core of $f$ is a disjoint union of two copies of $F$, the other is a round fold map as discussed in Theorem \ref{thm:7} (let $E$
 be a manifold diffeomorphic to the cylinder $S^{m-n} \times [-1,1]$ in the situation of Theorem \ref{thm:7}) and $M$ is a manifold as
 in the condition (\ref{thm:10.2}). 

 Conversely, from Theorem \ref{thm:1} (\ref{thm:1.1}), Example \ref{ex:2} (\ref{ex:2.1}) and Proposition \ref{prop:2}, a manifold as in the condition (\ref{thm:10.2}) always admits a round fold map satisfying the
 condition (\ref{thm:10.1}). 

 This completes the proof.
\end{proof}

\begin{Thm}
\label{thm:11}
Let $m, n>0$ be integers satisfying the relation $m \geq 2n$. For a manifold $M$ represented as a connected sum of $l>0$ manifolds regarded as the total spaces of smooth $S^{m-n}$-bundles over $S^n$ and $l+1$ {\rm (}n-1{\rm )}-connected
 manifolds admitting round fold maps into ${\mathbb{R}}^n$.
\end{Thm} 
\begin{proof}
We obtain a desired map by applying Corollary \ref{cor:1} starting from the pair of a map in Theorem \ref{thm:10} and a round fold map on an {\rm (}n-1{\rm )}-connected
 manifold in the assumption.
\end{proof}

\begin{Ex}
\label{ex:5}
A manifold $M$ represented as a connected sum of the following two manifolds admits a round fold map into
 the plane in
 Theorem \ref{thm:10}.
\begin{enumerate}
\item The total space of a smooth $S^4$-bundle over $S^2$.
\item The total space of a smooth ($S^2 \times S^2$)-bundle over $S^2$ presented in Example \ref{ex:3} (\ref{ex:3.2}).
\end{enumerate}
Moreover, a manifold represented as a connected sum of $l>0$ manifolds regarded as the total spaces of smooth $S^4$-bundles over $S^2$ and $l+1$ manifolds regarded as the total spaces of smooth ($S^2 \times S^2$)-bundles over $S^2$ admits a round fold map
 presented in Theorem \ref{thm:11}: note that the fibers of the latter $l+1$ manifolds (total spaces) may be arbitrary $4$-dimensional closed and simply-connected manifolds.
\end{Ex}
\begin{Rem}
There are infinitely many $6$-dimensional closed and simply-connected manifolds. For this, see also \cite{wall} etc..  
We don't know whether $6$-dimensional closed manifolds of wider classes, for example, $6$-dimensional closed and simply-connected manifolds admitting fold maps
 into the plane or equivalently having even Euler numbers, always admit round fold maps into the plane. 
\end{Rem}
The following is an extension of Theorem \ref{thm:9}.
\begin{Thm}
\label{thm:12}
Let $m,n>0$ be integers satisfying the relation $m \geq 2n$. Let $M$ be a closed and connected manifold of dimension $m$. Then
 the following two are equivalent.

\begin{enumerate}
\item
\label{thm:12.1}
 A round fold map $f:M \rightarrow {\mathbb{R}}^n$ satisfying the following two exist.
\begin{enumerate}
\item Inverse images of regular values of are disjoint unions of standard spheres and the number of connected components of the inverse images of a point in a proper core of $f$
 and the number of connected components of $S(f)$ are $l>1$ in the case $n \geq 1$ and $l>2$ in the case $n=1$. 
\item Let $f_0:M \rightarrow {\mathbb{R}}^n$ be a normal form of the round fold map $f$. The composition of the map 
 ${f_0} {\mid}_{{f_0}^{-1}({D^n}_{k+\frac{1}{2}}-{\rm Int} {D^n}_{k-\frac{1}{2}})}:{f_0}^{-1}({D^n}_{k+\frac{1}{2}}-{\rm Int} {D^n}_{k-\frac{1}{2}}) \rightarrow {D^n}_{k+\frac{1}{2}}-{\rm Int} {D^n}_{k-\frac{1}{2}}$
 and the map ${\pi}_{P}(x):=\frac{1}{2}\frac{x}{|x|}$ defines a bundle structure by Ehresmann's fibration theorem for
 any integer $1 \leq k \leq l-1$. The
 bundles are all trivial smooth bundles and the fibers are PL homeomorphic to disjoint unions of finite copies of the cylinder $S^{m-n} \times [-1,1]$ with
 the interior of a standard closed disc $D^{m-n+1}$ smoothly embedded in
 the interior of the manifold removed.
\item Consider the Reeb space $W_f$. For any connected component $R$ of $W_f-q_f(S(f))$ such that $\bar{f}(R)$ and every proper core of $f$ are disjoint, the inverse image ${q_f}^{-1}(a)$ of a point $a$ in $R$ is diffeomorphic to $S^{m-n}$.  
\end{enumerate}
\item
\label{thm:12.2}
$M$ is represented as a connected sum of $l-1$ manifolds regarded as the total spaces of smooth bundles whose fibers are almost-spheres of dimension $m-n$ over $S^n$. Further more, the number of the bundles whose fibers are not standard spheres are at least $l-1-k$ where $k$ is the maximal integer satisfying the relation $k \leq \frac{l}{2}$
\end{enumerate}
\end{Thm} 
\begin{proof}
We can show this in a way similar to the proof of Theorem \ref{thm:9}. We need to use maps in Example \ref{ex:2} (\ref{ex:2.1}) such that the almost-sphere may not be standard. For the number of the bundles in the assumption  (\ref{thm:12.2}), note that the number of connected components of the inverse image of a point in a proper core of the map $f$ is $l$.
\end{proof}

\end{document}